\documentclass{aims} 

\usepackage{amsthm}
\usepackage{amssymb}
\usepackage{mathtools}
\usepackage[noadjust]{cite}
\usepackage[misc]{ifsym}
\usepackage[shortlabels]{enumitem}
\usepackage{epsfig} 
\usepackage{epstopdf} 
\usepackage[colorlinks=true]{hyperref}
\usepackage[final]{changes}
\hypersetup{urlcolor=blue, citecolor=red}
\allowdisplaybreaks%

\def\currentvolume{X}
 \def\currentissue{X}
  \def\currentyear{200X}
   \def\currentmonth{XX}
    \def\ppages{X-XX}
     \def\DOI{10.3934/xx.xxxxxxx}

\newtheorem{theorem}{Theorem}[section]
\newtheorem{corollary}[theorem]{Corollary}

\newtheorem{lemma}[theorem]{Lemma}
\newtheorem{lemma-and-example}[theorem]{Lemma and Example}
\newtheorem{proposition}[theorem]{Proposition}

\theoremstyle{definition}
\newtheorem{definition}[theorem]{Definition}
\newtheorem{remark}[theorem]{Remark}
\newtheorem{assumption}[theorem]{Assumption}

\newtheorem{example}[theorem]{Example}

\DeclareMathOperator{\ad}{ad}
\DeclareMathOperator{\avg}{avg}
\DeclareMathOperator{\low}{\ensuremath\ell}
\DeclareMathOperator{\up}{\ensuremath\mathit{u}}

\newcommand{\R}{\mathbb{R}}  
\newcommand{\N}{\mathbb{N}}
\newcommand{\B}{\mathbb{B}}
\newcommand{\W}{\mathbb{W}}

\newcommand{\Lkal}{\mathcal{L}}

\newcommand{\Bkal}{\mathcal{B}}
\newcommand{\Akal}{\mathcal{A}}
\newcommand{\Fkal}{\mathcal{F}}
\newcommand{\Ekal}{\mathcal{E}}
\newcommand{\Ikal}{\mathcal{I}}
\newcommand{\Dkal}{\mathcal{D}}
\newcommand{\Mkal}{\mathcal{M}}

\DeclareMathOperator{\dist}{dist}
\DeclareMathOperator{\supp}{supp}
\DeclareMathOperator{\Dom}{Dom}
\DeclareMathOperator{\tr}{tr}
\renewcommand{\d}{\,\textnormal{d}}

\newcommand{\domLp}{\Dom(A+1)}
\newcommand{\irreg}{(IR)}
\newcommand{\reg}{(R)}
\newcommand{\altS}{\mathcal{S}_{\text{lin}}}
\newcommand{\Alin}{\Bkal}
\newcommand{\eps}{\varepsilon}

\DeclarePairedDelimiter{\norm}{\lVert}{\rVert}
\DeclarePairedDelimiter{\abs}{\lvert}{\rvert}

\numberwithin{equation}{section}

\title[Quasilinear parabolic control with gradient terms \& constraints]
{Optimal control of quasilinear parabolic PDEs with gradient terms and pointwise constraints on the gradient of the state}

\author[Lucas Bonifacius, Fabian Hoppe, Hannes Meinlschmidt and Ira Neitzel]{}

\subjclass{Primary: 49J20, 49K20, 35K59, 35K90; Secondary: 35B65, 35R05, 49K27}
\keywords{Optimal control, quasilinear parabolic PDE, gradient-state constraints, optimality conditions, existence of solutions}

\thanks{$^*$Corresponding author: Hannes Meinlschmidt}

\begin{document}
        \maketitle

        \centerline{\scshape
                Lucas Bonifacius $^{{\href{mailto:L.bonifacius@gmail.com}{\textrm{\Letter}}}}$,
                Fabian Hoppe $^{{\href{fabian.hoppe@dlr.com}{\textrm{\Letter}}}*1}$,
                Hannes Meinlschmidt $^{{\href{mailto:meinlschmidt@math.fau.de}{\textrm{\Letter}}}*2}$
                and Ira Neitzel $^{{\href{mailto:neitzel@ins.uni-bonn.de}{\textrm{\Letter}}}*3}$}

        \medskip

        {\footnotesize
                \centerline{$^1$German Aerospace Center (DLR), Institute of Software Technology,}\centerline{High-Performance Computing, Linder Höhe, Cologne, Germany}
        }

        \medskip

        {\footnotesize
                \centerline{$^2$Chair for Dynamics, Control, Machine Learning and Numerics (AvH professorship), Department of Mathematics}\centerline{Friedrich-Alexander-Universität Erlangen-Nürnberg, Cauerstraße 11, 91058 Erlangen, Germany}
        }

        \medskip
        {\footnotesize
                \centerline{$^3$Institut für Numerische Simulation, Rheinische Friedrich-Wilhelms-Universität Bonn,}\centerline{Friedrich-Hirze\-bruch-Allee 7, 53115 Bonn, Germany}
        }

        \bigskip



        \begin{abstract}
          We derive existence results and first order necessary optimality conditions for optimal control
          problems governed by quasilinear parabolic PDEs with a class of first order nonlinearities that
          include for instance quadratic gradient terms. Pointwise in space and time or averaged in space
          and pointwise in time constraints on the gradient of the state control the growth of the
          nonlinear terms. We rely on and extend the improved regularity analysis for quasilinear
          parabolic PDEs on a whole scale of function spaces from~\cite{HoppeMeinlschmidtNeitzel2023}.
          In case of integral in space gradient-constraints we derive first-order optimality conditions
          under rather general regularity assumptions for domain, coefficients, and boundary conditions,
          similar to~e.g.~\cite{Bonifacius2018}. In the case of pointwise in time and space
          gradient-constraints we use slightly stronger regularity assumptions leading to a classical
          smoother $W^{2,p}$-setting similar to~\cite{Casas2018}.
        \end{abstract}


\section{Introduction}
In the present paper we prove existence of solutions and first-order necessary optimality conditions for state-gradient-constrained optimal control problems of the form
        \begin{align} \tag{P}\label{ocp}
                \left\{ \qquad
                \begin{aligned}
                        &\min_{y,u} J(y,u) \coloneqq  \frac{1}{2}\lVert y - y_d \rVert_{L^2(I \times \Omega)}^2 + \frac{\gamma}{2} \lVert u \rVert_{L^2(\Lambda)}^2, \\[0.5em]
                        &\text{subject to} \quad u \in U_{\ad}, \quad y \in Y_{\ad},  \quad  (y,u)~\text{ solve }~\eqref{se},
                \end{aligned}
                \right.
        \end{align}
where~\eqref{se} is an abstract quasilinear parabolic PDE for the state $y$ and control $u$, of the general form
                \begin{align}\label{se}\tag{Eq}
                \partial_t y + \Akal(y) y = B u + \Fkal(y), \qquad
                y(0) = y_0.
        \end{align}

The abstract problem formulation also includes mixed boundary conditions and incorporates nonsmooth problem data.
The precise functional analytic setting, including the
definition of boundary conditions and the control-action operator $B$ will be introduced in Section~\ref{sec:notass}. The latter may in particular realize distributed controls (up to space dimension $d=3$), Neumann boundary controls (up to space dimension $d=2$), or purely time-dependent control functions coupled to fixed actuators with sufficient spatial regularity; see in particular some general settings in  Assumption~\ref{gaseq} and restrictions in Remark~\ref{rem-comments-on-irreg}.
Here, let us only highlight the key components and challenges of the optimization problem. First, \[\Akal(y)\coloneqq  - \nabla \cdot \xi(y) \mu \nabla\]
is a quasilinear second order elliptic operator and $\Fkal$ is a nonlinearity  of order at most one. Specific choices of $\Fkal$ may in particular include, but are not limited to,
\begin{equation}\label{eq:examplesF1}
\Fkal(y) = \lvert \nabla y \rvert_2^2,
\end{equation}
or
\begin{equation}\label{eq:examplesF2}
\Fkal(y) = y \beta \cdot \nabla y.\end{equation}

From an application point of view, the first example for $\Fkal$ arises in, e.g., the KPZ equation~\cite{KPZ1986} modelling the growth of interfaces, and it can be viewed as a prototype for certain nonlinearities arising in the modelling of semiconductors~\cite[Remark~4.7]{Meinlschmidt2021}. The second example can be considered as a prototype for the first-order nonlinear term from~\cite[Section~4]{Clement1992} related to the modelling of groundwater flows. Both are mathematically challenging because they grow quadratically in $y$ and thereby involve $\nabla y$, from which existence of solutions to the state equation becomes a delicate issue; see also~\cite[Section~6]{Dintelmann2009} for some discussion. In fact, there is a vast amount of blow-up results for $\Fkal$ of similar structure to~\eqref{eq:examplesF1} or~\eqref{eq:examplesF2}, for example for
\begin{equation}\label{eq:counterexF1}\Fkal(y) = \lambda y^r - \lvert \nabla y \rvert_2^2\end{equation}
with $r > 2$ and $\lambda > 0$ sufficiently large~\cite[Theorem~4]{Kawohl1989}, or
\begin{equation}\label{eq:counterexF2}\Fkal(y) = a \lvert \nabla y \rvert_2^{\alpha} + b y^{\beta}\end{equation}
with $a,b > 0$ and $\alpha,\beta > 1$~\cite[Theorem~3.2]{Hesaaraki2004}. A classical strategy to guarantee well-posedness of nonlinear PDEs is the use of growth conditions or monotonicity assumptions on the nonlinearity.
However, relying solely on growth conditions for $\Fkal$ limits the set of nonlinearities considerably. Alternatively, growth conditions can be combined with constraints on the state. This strategy has been used for instance in~\cite{Meinlschmidt2017_1,Meinlschmidt2017_2}, where pointwise bounds on the state combined with a tailored regularization in the objective function enforce global-in-time existence of a solution to the thermistor problem---a coupled system of a quasilinear parabolic and a quasilinear elliptic PDE\@.

We now take this idea a step further. We will allow nonlinearities $\Fkal$ that are uniformly bounded for all functions $y$ whose $W^{1,q}$-norm is uniformly bounded with respect to time, for some $q>d$; see Section~\ref{sec:notass} and in particular Assumption~\ref{gassemilin} for the precise setting. This is combined with certain gradient-state-constraints and box constraints on the control in Problem~\eqref{ocp} to control the growth of the nonlinear terms.
In essence, the constraints we consider are variations of either pointwise constraints on the gradient in space and time, i.e.
\begin{equation}\label{e:yadptwise}
    Y_{\ad} \coloneqq \Bigl\{ y\colon \lvert \nabla y \rvert_2^2 \leq g \ \text{ pointwise on } I \times \Omega \Bigr\},
\end{equation}
or averaged in space and pointwise in time gradient-constraints of the type
\begin{equation}\label{e:yadavergd}
    Y_{\ad} \coloneqq \Bigl\{ y\colon \lVert \nabla y \rVert^q_{L^q(\Omega)} \leq g_{\avg} \ \text{ pointwise on } I \Bigr\};
\end{equation}
with some possible generalizations commented on in Section~\ref{sec:generalize}.
Details will be provided in Assumption~\ref{gasconstr}, and a mathematically precise formulation of Problem~\eqref{ocp} will be stated in Section~\ref{sec::existence}.

The advantages and novel challenges of this approach are three-fold: First, gradient constraints are relevant in practical situations. For instance, when controlling the temperature of an object,
condition~\eqref{e:yadptwise} can be prescribed to guarantee that the Euclidean norm of the temperature gradient---which is related to thermal stresses in the material---everywhere in the object remains below a certain safety level at all times. Second, imposing gradient constraints in fact induces an implicit growth condition and allows to treat a considerably larger class of nonlinearities $\Fkal$ than relying solely on \emph{a priori} growth conditions, even when possibly combined with zero-order constraints. Third, however, first-order constraints on $\nabla y$ require higher regularity of the solutions of the PDE than zero-order constraints on $y$ as discussed in e.g.~\cite{Meinlschmidt2017_1,Meinlschmidt2017_2}. This is mainly due to classically relying on Slater-type arguments for the discussion of first order necessary conditions which in turn requires a function space setting where the set of admissible states has nonempty interior; cf.\ also the early work~\cite{Casas86} in the linear elliptic setting. For gradient-state constraints, we
will therefore need continuity of the gradient or the averaged in space gradient. The latter type relaxes the regularity requirements compared to fully pointwise constraints. This approach has been used in various settings in the literature.
Let us only mention~\cite{CasasFernandez1995} for a quasilinear elliptic problem with integrated constraints depending on the state and its gradient, and~\cite{Casas2001} for plain integral-gradient constraints for semilinear elliptic problems.
In the parabolic case,~\cite{Goldberg1993} have considered a semilinear problem with state constraints pointwise in time but integrated in space,
parabolic linear and semilinear problems with integrated in space and pointwise in time gradient bounds have been discussed in e.g.,~\cite{Mackenroth1983,Ludovici2015}.
Pointwise gradient constraints require more regularity. For the semilinear elliptic case we refer to~\cite{CasasFernandez1993}. Let us note in passing that there are several further contributions on finite element discretization and regularization methods for related
problems. Here, we only mention the discussion of a linear-quadratic parabolic problem with pointwise in time and integrated in space gradient constraints with purely time-dependent control functions in~\cite{Ludovici2015}, and do not give a detailed literature overview about numerical analysis results.

In the quasilinear setting of the present paper, the regularity of the state $y$ and potential blow-up-behavior is effected not only by semilinear terms $\Fkal$ but also the quasilinear operator $\Akal$. A suitable regularity setting also incorporating mixed boundary conditions and nonsmooth problem data without blow up is meanwhile well-established for quasilinear problems and even control problems with $\Fkal\equiv 0$ or nonlinearities $\Fkal$ that are usually monotone or grow at most linearly; see e.g.~\cite{Bonifacius2018,Casas2018,Meinlschmidt2016,HoppeMeinlschmidtNeitzel2023}. There, one considers the abstract equation~\eqref{se} in the negative Sobolev space $W^{-1,p}_D(\Omega)$, relying on optimal elliptic regularity in $W^{1,p}_D(\Omega)$ for the quasilinear elliptic operators. However, this setting does not suffice for our particular setting and the nonlinearities that inspired our research such as~\eqref{eq:examplesF1} and~\eqref{eq:examplesF2};
cf.\ the discussion in~\cite[Section 6]{Dintelmann2009}.
Instead, we essentially build upon existence and regularity results for quasilinear PDEs obtained in~\cite{HoppeMeinlschmidtNeitzel2023} in terms of nonautonomous maximal parabolic regularity for $\Akal$; cf.\ in particular Theorem~\ref{thm::existenceBoNe}, where $\Fkal$ did not depend on $y$ or $\nabla y$.

With the present work, we thus extend results for gradient-constrained optimal control problems to the quasilinear parabolic case. Building on the extensive regularity analysis of~\eqref{se} for $\Fkal \equiv 0$ on the scale of Bessel potential spaces $H^{-\zeta,p}_D(\Omega)$ with $\zeta \in [0,1]$ from~\cite{HoppeMeinlschmidtNeitzel2023} we can easily trade regularity assumptions on the equation against strength of the imposed gradient-constraints without duplicating arguments in a flexible manner: on the one hand, we are able to keep relatively low regularity assumptions on the state equation as in, e.g.,~\cite{Meinlschmidt2016,Bonifacius2018}, when considering integral in space and pointwise in time constraints on the gradient of the state. This setting, however, excludes the discussion of first-order conditions for pointwise in space and time gradient-constraints due to insufficient regularity of the states. On the other hand, we therefore also deal with a smoother setting for the state equation similar to~\cite{Casas2018}, and address first-order necessary optimality conditions for pointwise gradient-constraints.

From another point of view, we also extend the first-order theory on PDE-constrained optimization of quasilinear parabolic PDEs, where a lot of progress has been made in recent years with respect to the analysis of existence and regularity of solutions of the highly nonlinear state equation and its linearizations that is a typical difficulty in optimal control. We only mention some recent papers and refer the reader to the introduction of~\cite{Bonifacius2018} for an overview on quasilinear elliptic and earlier results on quasilinear parabolic problems. Wellposedness of a quasilinear parabolic state equation and existence of optimal controls for a respective control problem has been proven in~\cite{Meinlschmidt2016}
for nonsmooth domains and coefficients and mixed boundary conditions.
In~\cite{Bonifacius2018}, first- and second-order optimality conditions
have been derived for~\eqref{se} with $\Fkal \equiv 0$ subject to control constraints.
In~\cite{Casas2018}, a problem similiar to the one in~\cite{Bonifacius2018} is discussed: With respect to domain and coefficients the authors consider a smoother setting, but their nonlinearities are possibly unbounded and a monotone semilinearity $\Fkal$ of order zero is included. The results of~\cite{Bonifacius2018,Casas2018} have been extended to the presence of additional constraints on the state-variable, both pointwise in space and time and averaged in time but pointwise in space in~\cite{HoppeNeitzel2020}, for quasilinear PDEs with $\Fkal \equiv 0$. As mentioned before, regularity of solutions to~\eqref{se} with zero-order semilinearity on a whole scale of function spaces
has recently been investigated in~\cite{HoppeMeinlschmidtNeitzel2023}, connecting in some sense the regularity results from~\cite{Meinlschmidt2016,Bonifacius2018,Casas2018}; our analysis heavily builds on these results.
A central tool in the investigation of the appearing PDEs in this context is the functional analytic concept of maximal parabolic regularity for nonautonomous operators; cf.~\cite{Amann2004,Meinlschmidt2016}.

Our general approach is the following: A classical result from~\cite{Amann2005_3} will provide existence of local-in-time solutions of~\eqref{se}. The more recent results from~\cite{HoppeMeinlschmidtNeitzel2023} on nonautonomous maximal parabolic regularity guarantee that the quasilinear operator $\Akal$ does not cause any blow up. Combining these two results, we show that our assumptions on $\Fkal$ prevent blow up due to the semilinear term on a certain set of controls. Provided that this set is not empty, the interplay between the constraints prescribed in $Y_{\ad}$ and the properties of $\Fkal$ eventually guarantee well-posedness of the optimization problem and existence of a control-to-state operator with all necessary regularity properties. Important implications of these results are the following: First, for a given nonlinearity $\Fkal$, they allow to choose a set of constraints that guarantee well-posedness. Second, if gradient constraints are given that imply boundedness of the state in a $W^{1,q}$-norm for \textit{any} $q>d$, then our results provide an implicit growth condition on $\Fkal$ and thus a characterization of admissible nonlinearities, along with appropriate regularity assumptions. This is due to a unified theory for the parabolic PDEs on a scale of function spaces.

The rest of the paper is structured as follows: Section~\ref{sec:notass} introduces notation, definitions, and basic assumptions on the problem data. Auxiliary results are collected in Section~\ref{sec::auxresults} for later reference, while the quite technical verification of our assumptions for the model nonlinearities mentioned in the introduction has been moved to Appendix~\ref{sec:model-nonlinearity}. Section~\ref{sec::solutions} addresses solutions of the state equation~\eqref{se}: existence and uniqueness of local-in-time solutions is obtained in Section~\ref{sec::local-in-time}, and properties of global-in-time solutions (if they exist) are investigated in Section~\ref{sec:global-in-time}; up to this point our analysis exclusively deals with the state equation~\eqref{se}. Finally, from Section~\ref{sec::existence} on we turn to the control problem: we precisely introduce the gradient-constraints under consideration, formulate a corresponding additional assumption on the regularity setting, and prove existence of globally optimal solutions. Section~\ref{sec::fon} is then devoted to the derivation of first-order necessary optimality conditions, both for integral in space and pointwise in time gradient-constraints (Section~\ref{sec:fon-integr}) and pointwise in space and time constraints on the gradient (Section~\ref{sec::higher}) under the corresponding assumptions. Afterwards, some extensions and variations of our results are summarized in Section~\ref{sec:generalize}. For the convenience of the reader, frequently used definitions and results concerning the concept of maximal parabolic regularity have been collected in Appendix~\ref{sec:mpr}, and we verify that the model nonlinearities~\eqref{eq:examplesF1} and~\eqref{eq:examplesF2} satisfy the given array of assumptions in Appendix~\ref{sec:model-nonlinearity}.

\section{Notation, definitions, and assumptions}\label{sec:notass}

The goal of this section is to summarize the overall functional analytic setting and the assumptions on the problem data used in this paper. Moreover, we give concrete examples for admissible constellations.

\subsection{General Notation}

Given two Banach spaces $X$ and $Y$ we denote by $\Lkal(X,Y)$ the space of bounded linear maps $X \rightarrow Y$, equipped with the operator norm, and by $\Dom_X(A)$ the domain of a closed operator $A\colon X \rightarrow Y$, equipped with the graph norm. The space $X^* \coloneqq \Lkal(X,\R)$ stands for the topological dual of $X$. By $X \hookrightarrow Y$ we denote that $X$ embeds continuously into $Y$, and subscripts $d$ and $c$ indicate that this embedding is in addition dense or compact, respectively. Also, we use standard notation for the real $(X,Y)_{\theta,p}$ and complex interpolation spaces $[X,Y]_\theta$ between $X$ and $Y$. We refer to~\cite[Chapter~I.2]{Amann1995} for a concise resource. Given $p \in [1,\infty]$, the conjugate exponent of $p$ is denoted by $p'$, $1/p' + 1/p = 1$, with the usual conventions.

We will also heavily lean on the concept of maximal parabolic regularity. In order to not inflate this section too much with abstract results, we collect the necessary concepts and definitions in the appendix and only mention here that, for an interval $J$ and Banach spaces $Y \hookrightarrow_d X$, we set
\[
\W^{1,r}\bigl(J, X,Y\bigr) \coloneqq W^{1,r}(J, X) \cap L^r(J, Y).
\]

\subsection{Domain and function spaces}
For all of the following, the time interval $I = (0,T)$ with $T > 0$ is fixed, and we abbreviate the space time cylinder by $Q \coloneqq I \times \Omega$. The assumptions on $\Omega$ are as follows:
\begin{assumption}[Geometry]\label{ass:geometry}Let $\Omega \subset \R^d$, $d \in \{2,3\}$, be a bounded domain. Let further $\Gamma_N \subset \partial \Omega$ be a relatively open subset of $\partial\Omega$, the
\textbf{N}eumann boundary part, whereas $\Gamma_D \coloneqq \partial \Omega \setminus \Gamma_N$ denotes the \textbf{D}irichlet boundary part.
We assume that $\Omega \cup \Gamma_N$ is regular in the sense of Gröger~\cite{Groeger1989}, that is, a Lipschitz manifold (weak Lipschitz domain) with a compatiblity condition at the interface of $\Gamma_N$ and $\Gamma_D$, with the additional property that every boundary chart can be chosen to be volume-preserving.
\end{assumption}

We use standard notation for spaces of continuous, (continuously) differentiable, and Hölder continuous functions as well as for Lebesgue spaces. For $q \in (1,\infty)$ and $k \in \N$, let $W^{k,q}(\Omega)$ be the usual intrisically defined Sobolev spaces of order $k$, i.e., consisting of all functions $L^q(\Omega)$ whose distributional derivatives up to order $k$ are also in $L^q(\Omega)$, equipped with the natural standard sum norm. In addition, we define the integer Sobolev space $W^{1,q}_D(\Omega)$ incorporating a generalized zero trace condition on $D$ by taking the closure of $C_D^\infty(\Omega)$ in $W^{1,q}(\Omega)$, where
\begin{equation*}
  C_D^\infty(\R^d) \coloneqq \Bigl\{ \varphi \in C_c^\infty(\R^d),~\dist(\supp \varphi,\Gamma_D) > 0\Bigr\}, \qquad C_D^\infty(\Omega) \coloneqq C_D^\infty(\R^d)\bigr|_\Omega.
\end{equation*}
In order to also have Sobolev type spaces with noninteger differentiability at our disposal, we consider the scale of Bessel potential spaces on $\Omega$, also incorporating a generalized zero trace condition on $D$, as follows: For $s \in \R$ and $q \in (1,\infty)$, let $H^{s,q}(\R^d)$ be the usual scale of Bessel potential spaces on $\R^d$. We take these for granted.
Then from~\cite[Theorem~VII.1]{Jonsson1984}, for $s > 1/q$, there exists a continuous linear trace operator $\tr_D \colon H^{s,q}(\R^d) \to L^q(\Gamma_D)$, the latter with respect to the $(d-1)$-dimensional Hausdorff measure. Let $s \geq 0$. Define $H^{s,q}_D(\R^d) \coloneqq \ker\tr_D$ if $s>1/q$ and $H^{s,q}_D(\R^d) \coloneqq H^{s,q}(\R^d)$ otherwise and finally set
\begin{equation*}
  H^{s,q}(\Omega) \coloneqq H^{s,q}(\R^d)\bigr|_\Omega, \qquad H^{s,q}_D(\Omega) \coloneqq H^{s,q}_D(\R^d)\bigr|_\Omega.
\end{equation*}
Note that for $s < 1/q$, by definition $H^{s,q}_D(\Omega) = H^{s,q}(\Omega)$, and that $H^{s,q}(\Omega) = L^q(\Omega)$ for $s=0$. With the canonical quotient norm, these spaces defined by restriction are complete. We further set $W^{-1,q}_D(\Omega) \coloneqq W^{1,q'}_D(\Omega)^*$ and $H^{-s,q}_D(\Omega) \coloneqq H^{s,q'}_D(\Omega)^*$. Further, the geometry assumption that $\Omega$ is a weak Lipschitz domain implies that in fact $H^{k,q}(\Omega) = W^{k,q}(\Omega)$ for integer $k$ and also $H^{1,q}_D(\Omega) = W^{1,q}_D(\Omega)$, each up to equivalent norms~(\cite[Proposition~B3]{Bechtel2019}), since there is a universal extension operator. As a consequence, $W^{1,q}_D(\Omega) \cap W^{1,r}(\Omega) = W^{1,r}_D(\Omega)$ for any $r \in (1,\infty)$. For $s \in (0,1)$, we also recover $H^{s,q}_D(\Omega)$ and $H^{s,q}(\Omega)$ by interpolation from the familiar Lebesgue- and Sobolev spaces:
\begin{lemma}[{\cite[Theorem~1.1]{Bechtel2019}}]\label{lem:bessel-interpolation-scale-identification} Let $q \in (1,\infty)$. Then for every $\theta \in (0,1) \setminus \{\frac1q\}$, we have, up to equivalent norms,
  \begin{equation*}
\Bigl[L^q(\Omega),W^{\pm 1,q}_D(\Omega)\Bigr]_\theta = H^{\pm\theta,q}_D(\Omega)
  \end{equation*}
  and
\begin{equation*}
\Bigl[L^q(\Omega),W^{1,q}(\Omega)\Bigr]_\theta = H^{\theta,q}(\Omega).
  \end{equation*}
\end{lemma}
Due to the definition of $H^{s,q}(\Omega)$ and $H^{s,q}_D(\Omega)$ as quotient spaces, we also have the usual Sobolev embeddings for these spaces at hand, be that within the Bessel scale or into spaces of $1 \wedge (s-d/q)$-H\"older continuous functions on $\overline\Omega$ if $s > d/q$.
Finally, let us note that since $\Omega$ is fixed throughout the paper, we will usually
omit it in the notation of function spaces and write $L^p$, $W^{1,p}_D$, $H^{\theta,p}$ etc.\ instead of $L^p(\Omega)$, $W^{1,p}_D(\Omega)$, $H^{\theta,p}(\Omega)$ etc.\ for the reason of brevity.

\subsection{The elliptic differential operator}

We next define the fundamental elliptic differential operator in divergence form $-\nabla \cdot\mu\nabla$. Nonlinear multiplicative perturbations in the coefficient function of this operator will give rise to the quasilinear versions defined thereafter; cf.~\eqref{eq:lin-A-def}. The basic assumption on the coefficient function $\mu$, valid for all of the following, is as follows.
\begin{assumption}\label{ass:coeff-matrix} Let $\mu\colon \Omega \rightarrow \R^{d\times d}$
be a measurable, uniformly bounded and coercive coefficient matrix function, in the following sense:
\[ \mu_{\bullet} \coloneqq \inf_{x \in \Omega} \inf_{z \in \R^d \setminus \{0\}} \frac{z^T \mu(x) z}{z^T z} > 0, \qquad \mu^{\bullet} \coloneqq \sup_{x \in \Omega} \sup_{1 \leq i,j \leq d} \lvert \mu_{i,j}(x) \rvert < \infty. \]
\end{assumption}
Note that we do not assume any smoothness of the coefficient function $\mu$.

We define $-\nabla \cdot \mu \nabla \colon W^{1,2}_D(\Omega) \to W^{-1,2}_D(\Omega)$ by
\begin{equation*}
\bigl\langle -\nabla\cdot\mu\nabla \varphi,\psi\bigr\rangle \coloneqq \int_\Omega \mu\nabla \varphi \cdot \nabla \psi \qquad (\varphi,\psi \in W^{1,2}_D(\Omega)).
\end{equation*}
With H\"older's inequality, it is clear that $-\nabla \cdot \mu\nabla$ is bounded. Further, from coercivity of $\mu$, we see that $-\nabla \cdot \mu \nabla + 1 \colon W^{1,2}_D(\Omega) \to W^{-1,2}_D(\Omega)$ is continuously invertible by the Lax-Milgram lemma. (Note that one can omit the ``$+1$'' if $\Gamma_D \neq \emptyset$ in the setting of Assumption~\ref{ass:geometry}.) For $q \neq 2$, we further consider the part of $-\nabla \cdot \mu \nabla$ in $W^{-1,q}_D(\Omega)$. In general, the domain of that operator is larger than $W^{1,q}_D(\Omega)$. Thus, we pose the following optimal elliptic regularity assumption, tailored to the particular setting of this work:
\begin{assumption}\label{ass:optimal-elliptic-regularity-mu}
We assume that there is $p >d$ such that
\begin{align*}
-\nabla \cdot \mu \nabla + 1\colon W^{1,p}_D(\Omega) \rightarrow W^{-1,p}_D(\Omega) \quad \text{is a topological isomorphism.}
\end{align*}
We fix this choice of $p$ now and for the rest of this paper. Consequently, we also choose and fix a number $s \in (1,\infty)$ such that $\frac2s < 1-\frac{d}p$.
\end{assumption}
Note that the foregoing assumption is nontrivial; it is an implicit restriction on $\mu$ versus the
geometry of $\Omega$, respectively of $\Gamma_D$ and $\Gamma_N$. Still, in the context of
Assumption~\ref{ass:geometry}, it is always true for $d=2 $ as proven in the seminal work of
Gröger~\cite{Groeger1989}, thus, it is in fact only an assumption for $d \geq 3$. See~\cite{Disser2015}
for a comprehensive study for $d=3$. Assumption~\ref{ass:optimal-elliptic-regularity-mu} will provide us
with stability of domains of the quasilinear counterparts of $-\nabla\cdot\mu\nabla + 1$, however, there
are many more useful consequences of which we collect a few. We are going to be quite concise here
and refer to~\cite[Section~2]{HoppeMeinlschmidtNeitzel2023} for a more detailed collection of results and
further comments, also on Assumption~\ref{ass:optimal-elliptic-regularity-mu} in general, as well as
references to other works. However, in Examples~\ref{ex:data-regularity-strong}
and~\ref{ex:data-regularity2} at the end of this section, we provide two exemplary settings that fulfill
the foregoing and following assumptions.

Denote by $A$ the part of $-\nabla \cdot \mu \nabla$ in $L^p(\Omega)$. Then both $A+1$ and $-\nabla \cdot \mu\nabla+1$ are \emph{positive} operators on $L^p(\Omega)$ and $W^{-1,p}_D(\Omega)$, respectively; in particular, their fractional powers are well defined. In fact, the operators even admit \emph{bounded imaginary powers} and thus satisfy \emph{maximal parabolic regularity}. It follows that also the part of $-\nabla \cdot \mu\nabla + 1$ in any of the interpolation spaces $[L^p(\Omega),W^{-1,p}_D(\Omega)]_\zeta$ for $\zeta \in [0,1]$ satisfies maximal parabolic regularity~\cite[Theorem~5.16]{Dintelmann2009}. Recall that the latter spaces coincide with $H^{-\zeta,p}_D(\Omega)$ by Lemma~\ref{lem:bessel-interpolation-scale-identification}. We will extensively use the abbreviation
\[\Dkal_\zeta \coloneqq \Dom_{H^{-\zeta,p}_D}(-\nabla \cdot \mu \nabla+1)\] for the domain of the part of $-\nabla \cdot \mu \nabla+1$ in $H^{-\zeta,p}_D$.

For the most important fractional power of $A+1$, we have the following fundamental result under Assumption~\ref{ass:optimal-elliptic-regularity-mu}, the \emph{Kato square root property}~(\cite{Bechtel2024},~\cite[Theorem~6.5]{terElst2017}):
\begin{proposition}\label{prop:kato} We have $\Dom (A+1)^{1/2} = W^{1,p}_D(\Omega)$ with equivalent norms, i.e.,
\begin{equation}
 \label{eq:kato-general} (A+1)^{1/2} \colon W^{1,p}_D(\Omega) \to L^p(\Omega) \quad \text{is a topological isomorphism.}
\end{equation}
\end{proposition}
The Kato square root property~\eqref{eq:kato-general} and several techniques that make use of it will prove invaluable over the course of this work; the attentive reader is encouraged to look for it in Section~\ref{sec::auxresults}. We give one example here:
\begin{lemma}\label{lem:regular-case-interpolation}
For $0 \leq \alpha \leq 1$, we have, up to equivalent norms,
\begin{equation}\label{eq:fracpower-equals-bessel-below-12}
\domLp^{\alpha/2}= H^{\alpha,p}_D.
\end{equation}
\end{lemma}

\begin{proof}
  We use the complex interpolation of fractional power domains of operators admitting bounded imaginary powers~(\cite[I.(2.9.8)]{Amann1995}). With the square root property $\domLp^{1/2} = W^{1,p}_D$ as in Proposition~\ref{prop:kato} and Lemma~\ref{lem:bessel-interpolation-scale-identification}, we get
  \begin{equation*}
\domLp^{\alpha/2}= \Bigl[L^p,\domLp^{1/2}\Bigr]_{\alpha} = \bigl[L^p,W^{1,p}_D\bigr]_{\alpha} = H^{\alpha,p}_D,
  \end{equation*}
  as claimed.
\end{proof}

Finally, we set the stage for the quasilinear elliptic operator in~\eqref{se}, by introducing the nonlinear function $\xi$.
\begin{assumption}\label{ass:xi}
Let $\xi\colon \R \rightarrow \R$ be a continuously differentiable real function with $0 < \xi_{\bullet} \leq \xi(x) \leq \xi^{\bullet}$ for all $x \in \R$. Further assume that $\xi'$ is Lipschitz continuous.
\end{assumption}

\begin{remark}
 Assumption~\ref{ass:xi} is quite standard and will allow us to have $y \mapsto \xi(y)$ and $y \mapsto \xi'(y)$ map H\"older and Lipschitz spaces into themselves, for which we cannot relax it much more. At the same time, we consider it to be not particularly restrictive in applications, where $\xi$ is used to model, say, a temperature-dependent thermal conductivity of the materials occupying the region $\Omega$. For example, for stainless steel, the thermal conductivity $\xi(y)$ can be represented by an affine linear function of $y$~(\cite[Section~4]{CRST93}) which we then smoothly cut off outside of physically relevant regions to make it formally satisfy Assumption~\ref{ass:xi}.
\end{remark}

With Assumption~\ref{ass:xi} at hand, we define the operators $\Akal$ and $\Alin$ in an \emph{ad hoc} manner by setting
\begin{equation}\label{eq:lin-A-def}
 \Akal(y) \coloneqq -\nabla \cdot \xi(y)\mu \nabla, \qquad \Alin(y) \coloneqq \Bigl[\varphi \mapsto -\nabla \cdot \xi'(y)\varphi\mu\nabla y\Bigr]
\end{equation}
in analogy to the definition of $-\nabla\cdot\mu\nabla$ \emph{mutatis mutandis}. The main part of the differential operator in~\eqref{se} is then precisely $\Akal(y)y$, and the linearization of this operator is given by $\Akal(y) + \Alin(y)$, which explains the significance of $\Alin$. Here, $y$ can (and will) be also a function involving the time variable, say, $I \times \Omega \to \R$, and we will tacitly consider $\Akal(y)$ and $\Alin(y)$ as operators acting on time-dependent abstract functions with the usual convention $(\Akal(y)z)(t) \coloneqq \Akal(y(t))z(t)$.

The crucial functional analytic property of $\Akal$ that we rely on later will be the following, which is in fact one of the main results in~\cite{HoppeMeinlschmidtNeitzel2023}; see Appendix~\ref{sec:mpr} for the notion of maximal regularity.
\begin{proposition}\label{prop:A-mpr-first-statement} Suppose that $\zeta\in[0,1]$ and $y \in C(\overline I,C^\sigma(\overline\Omega))$ and $\sigma>1-\zeta$ with $\sigma=1$ if $\zeta =0$. Then $\Akal(y)$ satisfies nonautonomous maximal parabolic regularity on $L^r(I,H^{-\zeta,p}_D(\Omega))$ with constant domain $\Dkal_\zeta$ for every $r \in (1,\infty)$.
\end{proposition}
We refer to~\cite[Lemma~4.2/Theorem~4.3]{HoppeMeinlschmidtNeitzel2023} for a proof; see also Lemma~\ref{lem::mprhminuszeta} below, where we also
establish nonautonomous maximal parabolic regularity for $\Akal(y) + \Alin(y)$ and comment on the discrepancy between the definition of $\Akal(y)$ and the ``$+1$'' in $\Dkal_\zeta$. These results will in fact be the workhorse that allows us to deal with the quasilinear equation~\eqref{se} in a flexible and convenient manner in a scale of function spaces. The crucial information is the uniformity of the domain of $\Akal(y)+1$ in $H^{-\zeta,p}_D$, which will be $\Dkal_\zeta$, with respect to $y$ in a suitable regularity class $\subseteq C(\overline I,C^\sigma(\overline\Omega))$. This allows us to concentrate on $\Dkal_\zeta$, and indeed, to set the stage for Lemma~\ref{lem::mprhminuszeta} and the optimal control considerations in the later sections, we invest some work in the domains $\Dkal_\zeta$ and associated function spaces in Section~\ref{sec::auxresults}.

\subsection{Regularity setting}

We consider two different regularity settings in this work, one being the \emph{irregular} one~\irreg, which was introduced so far, but for which we have to restrict $\zeta > \frac{d}p$ for several results due to structural reasons, and a \emph{regular} one~\reg. For the irregular setting, we can deal with the optimal control problem subject to pointwise-in-time but spatially integrated gradient constraints, whereas the regular one allows to deal with actual pointwise constraints on the gradient in time and space. However, the model nonlinearity $\abs{\nabla y}^2_2$ as in~\eqref{eq:examplesF1} in the introduction, and also~\eqref{eq:examplesF2}, can be dealt with in both regularity settings.
\begin{assumption}[Range of $\zeta$, regularity setting]\label{assu:regularity-settings}
Let $\zeta \in [0,1]$ and suppose at least one of the following additional conditions:
\begin{enumerate}
  \item[\irreg] $\zeta > \frac{d}p$, \emph{or}
  \item[\reg] $\domLp = W^{2,p}(\Omega) \cap W^{1,p}_D(\Omega)$.
\end{enumerate}
\end{assumption}

For a discussion of Assumptions~\ref{ass:geometry},~\ref{ass:coeff-matrix},~\ref{ass:optimal-elliptic-regularity-mu}, and~\ref{assu:regularity-settings} in the case \irreg, including possible non-trivial and non-smooth configurations that satisfy them, we refer to~\cite[Section 2.2]{Bonifacius2018}. The regular case \reg{} essentially implicitly restricts us to quite strong problem data, for example, to the following classical regular setting:

\begin{example}\label{ex:data-regularity-strong}
Assume that $\Omega$ has a $C^{1,1}$-boundary, $\mu$ satisfying Assumption~\ref{ass:coeff-matrix} is additionally Lipschitz continuous, and $\overline{\Gamma_D} \mathrel\cap \overline{\Gamma_N} = \emptyset$.\footnote{This does not necessarily imply that $\Gamma_D = \partial \Omega$ or $\Gamma_D = \emptyset$;
think, for example, of the annulus $\Omega = B_1 \setminus \overline{B_{1/2}}$ with $\Gamma_D = \partial B_1$ and $\Gamma_N = \partial B_{1/2}$ which is indeed admissible for Assumption~\ref{assu:regularity-settings} \reg{} even though two different types of boundary conditions are implied.} Then, Assumptions~\ref{ass:geometry},~\ref{ass:optimal-elliptic-regularity-mu}, and~\ref{assu:regularity-settings} \reg{} are fulfilled; cf.~\cite[Theorem~3.12, Remark~3.17]{Elschner2007} and~\cite[Theorems~2.4.2.5 and~2.4.2.7]{Grisvard1985}.
In space dimension $d=2$ this is also true for convex domains (without other assumptions on boundary smoothness) if $\Gamma_D = \partial \Omega$ or $\Gamma_D = \emptyset$; cf., e.g.,~\cite[Example~5.2]{HoppeMeinlschmidtNeitzel2023}.
\end{example}
The following example is an exemplary situation of case \irreg.
\begin{example}\label{ex:data-regularity2} Let $d=3$ and let $\Omega = (-1,1)^2 \times (-2,2)$ be a box. Choose the Neumann boundary part $\Gamma_N = \bigl((-\frac12,\frac12)^2 \times \{2\}\bigr) \cup \bigl((-\frac12,\frac12)^2 \times \{-2\}\bigr)$ to be a square on the bottom and top boundary parts of $\Omega$, and put $\Gamma_D = \partial\Omega\setminus \Gamma_N$. Consider a coefficient function $\mu = \chi_{B_1} + 2 \chi_{\Omega\setminus B_1}$ with the unit ball $B_1 \subseteq \R^3$. Note that $\mu$ has a jump discontinuity along $\partial B_1$ which touches $\Gamma_D \subset \partial\Omega$. Yet Assumptions~\ref{ass:geometry} and~\ref{ass:optimal-elliptic-regularity-mu} are satisfied, see~\cite[Theorem~4.8]{Disser2015} for the latter, but~\reg{} will not hold.
\end{example}

We chose an example where the jump discontinuity of $\mu$ does not interact with mixed boundary conditions in order to keep the technical assumptions to a reasonable level, but such cases are also admissible for Assumption~\ref{ass:optimal-elliptic-regularity-mu} in particular situations. Again, we refer to~\cite[Theorem~4.8]{Disser2015}.
\subsection{Problem data}
We collect the assumptions on the remaining problem data in~\eqref{ocp}, starting with the objective function and admissible set controls as well as initial condition for the state equation that remains fixed. For this, we fix here that Assumptions~\ref{ass:geometry},~\ref{ass:coeff-matrix},~\ref{ass:optimal-elliptic-regularity-mu},~\ref{ass:xi},~\ref{assu:regularity-settings} and the upcoming Assumption~\ref{gaseq} remain valid throughout \emph{all that follows}. In particular, we recall the following fundamental relations that will be used consistently in the following:
\[p>d,\qquad \frac{2}{s}<1-\frac{d}{p},\]
and $\zeta\in[0,1]$ fulfilling Assumption~\ref{assu:regularity-settings}. In contrast, for Assumption~\ref{gassemilin} below, we will activate assumptions on an individual basis.
\begin{assumption}[General problem data]\label{gaseq} For general data in~\eqref{ocp}, we collect:
\begin{enumerate}[(i)]
\item In~\eqref{se}, we consider an initial condition $y_0 \in (H^{-\zeta,p}_D(\Omega), \Dkal_\zeta)_{1-\frac1s,s}$.
\item For a finite measure space $(\Lambda, \rho)$ we denote the control space by $L^s(\Lambda)$ and define the set of admissible controls as
\[U_{\ad} = \Bigl\{ u \in L^s(\Lambda)\colon u_{\low}(x) \leq u(x) \leq u_{\up}(x) \quad \text{for $\rho$-a.a.}~x \in \Lambda\Bigr\}\]
with $u_{\low}, u_{\up} \in L^\infty(\Lambda)$, $u_{\low} \leq u_{\up}$ almost everywhere. The control action operator
\begin{align*}
B\colon L^s(\Lambda) \rightarrow L^s(I,H^{-\zeta,p}_D)
\end{align*}
is linear and bounded.
\item A desired state $y_d \in L^\infty(I,L^2(\Omega))$ and the cost parameter $\gamma > 0$ are fixed.
\end{enumerate}
\end{assumption}
We note that this setting allows to incorporate distributed controls (for each $\zeta$), purely time-dependent controls with actuators in $H^{-\zeta,p}_D$ (for each $\zeta$), as well as Neumann boundary control (for $\zeta > \frac{d}{p}$ and $d=2$ only); cf.\ the discussion in~\cite{Bonifacius2018}.

Next, let us collect the assumptions on the nonlinearity $\Fkal$ in~\eqref{se}.
\begin{assumption}[Nonlinear function $\Fkal$]\label{gassemilin}
We assume that $\Fkal \colon \W^{1,s}(I,H^{-\zeta,p}_D,\Dkal_\zeta) \to L^s(I,H^{-\zeta,p}_D)$ with the following additional properties:
\begin{enumerate}[(i)]
\item\label{gassemilin:volterra} It is a Volterra map such that $\Fkal - \Fkal(0)$ is Lipschitz continuous on bounded sets $\W^{1,s}(I,H^{-\zeta,p}_D,\Dkal_\zeta) \to L^r(I,H^{-\zeta,p}_D)$ for some $r>s$.

\item\label{gassemilin:boundedness} There is a $p \leq q \leq \infty$ such that
for every
$C\geq 0$ there is $M \geq 0$ such that for every $y \in \W^{1,s}(I,H^{-\zeta,p}_D,\Dkal_\zeta)$:
\begin{equation*}
  \sup_{t\in \overline I} \norm{y(t)}_{W^{1,q}_D} \leq C \quad \implies \quad \norm{\Fkal(y)}_{L^s(I,H^{-\zeta,p}_D)} \leq M.
\end{equation*}

In the irregular case~\irreg{} of Assumption~\ref{assu:regularity-settings}, we require the foregoing precisely for $q=p$.

\item\label{gassemilin:weak-cont} $\Fkal$ is weakly continuous.
\item\label{gassemilin:diffbar} $\Fkal$ is Fr\'{e}chet-differentiable and for every $y \in \W^{1,s}(I,H^{-\zeta,p}_D,\Dkal_\zeta)$, we have
\begin{equation*}
  \Fkal'(y) \in L^s\bigl(I,\Lkal((H^{-\zeta,p}_D,\Dkal_\zeta)_{1-\frac1s,s},H^{-\zeta,p}_D)\bigr).
\end{equation*}
\end{enumerate}
\end{assumption}

Before discussing these assumptions on $\Fkal$ below, let us immediately note that the two prototypical examples~\eqref{eq:examplesF1} and~\eqref{eq:examplesF2} for $\Fkal$,  so
$\Fkal(y) = \lvert \nabla y \vert_2^2$ and $\Fkal(y) = y\beta \cdot \nabla y$, indeed satisfy Assumption~\ref{gassemilin}, both also already in the irregular case~\irreg{} and with $q=p$ in Assumption~\ref{gassemilin}~\ref{gassemilin:boundedness}. The corresponding details which verify this can be found in Appendix~\ref{sec:model-nonlinearity}.

We point out that part~\ref{gassemilin:volterra} of Assumption~\ref{gassemilin} is needed in \emph{all that follows}, whereas the remaining Assumptions~\ref{gassemilin}~\ref{gassemilin:boundedness},~\ref{gassemilin:weak-cont} and~\ref{gassemilin:diffbar} are only needed in the discussion of the optimal control problem~\eqref{ocp}, starting from Section~\ref{sec:global-control-prop}. To be even more precise, the differentiability assumption in Assumption~\ref{gassemilin}~\ref{gassemilin:diffbar} is only required for the dicussion of first order necessary conditions in Section~\ref{sec::fon}.
\begin{remark}\label{rem:assumptions-F} A few comments on Assumption~\ref{gassemilin} are in order.
\begin{enumerate}[(i)]
  \item We do not focus on nonlocal-in-time equations in the present work, so we include the Volterra property~\ref{gassemilin:volterra} for formal reasons to use the main result in~\cite{Amann2005_3} later in our Proposition~\ref{prop::localex}. We refer to~\cite{Amann2005_3,Amann2005_4} for more details.
  \item To see the importance of Assumption~\ref{gassemilin}~\ref{gassemilin:diffbar}, we would like to point the reader to Proposition~\ref{prop:cont-mpr-plus-pruess} in the appendix.
  \item The growth condition in~\ref{gassemilin:boundedness} is quite weak and will be satisfied, for sufficiently large $q$, by nearly every nonlinear function acting in $y$ or $\nabla y$. In the irregular case~\irreg{}, where we are restricted to $q=p$, there \emph{is} a certain restriction but it is still quite general. We comment on prototypical examples for $\Fkal$ in Appendix~\ref{sec:model-nonlinearity}.
  \item The number $q$ in~\ref{gassemilin:boundedness} is a quite critical quantity in this work, as we will pose gradient constraints in the context of the optimal control problem~\eqref{ocp} exactly such that they enforce $\norm{y(t)}_{W^{1,q}_D}$ to be uniformly bounded; see Assumption~\ref{gasconstr} below. It is precisely this entanglement which allows us to deal with strong nonlinearities in~\eqref{se} in~\eqref{ocp}. However, in order to obtain a satisfying optimal control theory, starting from Section~\ref{sec::existence}, we will indeed have to link $\zeta$ and $q$ by requiring that \begin{equation}\label{zetaqlink}\zeta < 1-\frac2s - \frac{d}p + \frac{d}q\end{equation}
  to the effect that $\W^{1,s}(I,H^{-\zeta,p}_D,\Dkal_\zeta) \hookrightarrow C(\overline I,W^{1,q}_D(\Omega))$. The precise assumption will be stated in the box at the beginning of Section~\ref{sec::existence}.
\end{enumerate}
\end{remark}

We proceed with the formal definition of the constraint set $Y_{\ad}$ in~\eqref{ocp}. Here, we consider two prototypical sets of admissible states $Y_{\ad}$ defined by gradient constraints. The crucial aspect of the following assumption is that we \emph{link} the function spaces in which we pose the gradient constraints to the growth bound in Assumption~\ref{gassemilin}~\ref{gassemilin:boundedness}. This will eventually allow us to obtain a well posed optimal control problem~\eqref{ocp} for which we can develop a satisfying theory. As mentioned in Remark~\ref{rem:assumptions-F}, we will eventually establish a link between $\zeta$ and $q$ in Sections~\ref{sec::existence} and following, which is where the following assumption on the precise constraint setting will be put into practice.
\begin{assumption}\label{gasconstr}
Let $q$ be as in Assumption~\ref{gassemilin}~\ref{gassemilin:boundedness}. We assume that $\Gamma_D\neq\emptyset$ and define the set of admissible states to be one of the following:
\begin{enumerate}[(i)]
\item\label{gasconstr-avg} Integral in space and pointwise in time gradient constraints:
\begin{equation}\label{Yadav} Y_{\ad} \coloneqq \left\{ y \in \W^{1,s}(I,H^{-\zeta,p}_D, \Dkal_\zeta)\colon\lVert \nabla y \rVert_{L^q}^q \leq g_{\avg} \text{ a.e.\ on } I \right\} \end{equation}
with $g_{\avg} \in L^\infty(I)$ and $g_{\avg} \geq 0$ a.e.\ on $I$.

\item\label{gasconstr-ptw} Pointwise in space and time gradient constraints:
\begin{equation}\label{Yadptw} Y_{\ad} \coloneqq \left\{ y \in \W^{1,s}(I,H^{-\zeta,p}_D, \Dkal_\zeta)\colon \lvert \nabla y \rvert_2^2 \leq g \text{ a.e.\ on } Q \right\} \end{equation}
with $g \in L^\infty(I,L^{q/2})$ with $g \geq 0$ a.e.\ on $Q$.

\end{enumerate}
If $q = \infty$, we tacitly identify~\eqref{Yadav} with~~\eqref{Yadptw}.
\end{assumption}
To prevent any confusion or ambiguity, let us clarify that the $L^q$-norm of the (weak) gradient of a function used in~\eqref{Yadav} is defined as follows:
\begin{align}\label{eq:grad-lq-norm}
\lVert \nabla y \rVert_{L^q}^q := \int_\Omega \lvert \nabla y(x) \rvert_2^q \d x := \int_\Omega \Bigl(\sum_{i=1}^d [\partial_{x_i} y](x)^2 \Bigr)^{\frac{q}{2}} \d x
\end{align}

Several modifications of these constraints can be considered. However, we concentrate on~\eqref{Yadav} and~\eqref{Yadptw} as prototypical examples and postpone the discuss of some further variants, including pure Neumann boundary conditions with $\Gamma_D=\emptyset$ to Section~\ref{sec:generalize}.

\begin{remark}\label{rem:semibounds}
We point out that for our current setting with $\Gamma_D\neq\emptyset$, the pure gradient bounds in the formulation of~\eqref{Yadptw} and~\eqref{Yadav} allow to obtain ``full'' $W^{1,q}$-bounds as in the formulation of Assumption~\ref{gassemilin}~\ref{gassemilin:boundedness} due to the $L^q$-Poincaré-Friedrichs inequality~\cite[Lemma 1.36]{Gajewski1974}.
\end{remark}

\section{Auxiliary results}\label{sec::auxresults}

 The section contains known or straight forward auxiliary results for proofs of our main results on
 existence of optimal controls in Sections~\ref{sec::existence} and on first order optimality conditions
 in Sections~\ref{sec:fon-integr} and~\ref{sec::higher}. We recall that Assumptions~\ref{ass:geometry},~\ref{ass:coeff-matrix},~\ref{ass:optimal-elliptic-regularity-mu}, and~\ref{ass:xi} are presumed without further mention.

We begin with a standard embedding for the maximal regularity spaces $\W^{1,s}(I,X,Y)$ that we will use repeatedly in the following. It is a consequence of~\cite[(5.8)]{Amann2004}, stated in a convenient manner:

\begin{lemma}\label{lem:max-reg-embeddings} Let $1 < r < \infty$ and $\theta < 1-\frac1r$ as well as $\alpha < 1-\frac1r-\theta$. Then
\begin{equation}
 \label{eq:maxreg-embed}
   \W^{1,r}(I,X,Y) \hookrightarrow C^\alpha\bigl(\overline I,[X,Y]_\theta\bigr)
\end{equation}
and the embedding is compact if $Y \hookrightarrow X$ compactly.
\end{lemma}

In order to identify the interpolation spaces occurring in~\eqref{eq:maxreg-embed} for the choices $X = H^{-\zeta,p}_D$ and $Y = \Dkal_\zeta$, we further recall the following characterization of complex interpolation spaces between these spaces in terms of fractional power domains of the operator $-\nabla \cdot\mu\nabla + 1$, as stated in~\cite[Lemma~5.3]{HoppeMeinlschmidtNeitzel2023}.

\begin{lemma}\label{lem:interpolation-identify-frac-power}
Let $\frac\zeta2 < \theta \leq 1$. Then
\begin{equation}
 \label{eq:interpolation-identify-frac-power} \bigl[H^{-\zeta,p}_D,\Dkal_\zeta\bigr]_\theta = \domLp^{\theta-\zeta/2}.
\end{equation}
\end{lemma}

The foregoing lemmata are put to use in the next embedding result which is used several times throughout this work in order to deal with gradient nonlinearites and gradient constraints---see Corollaries~\ref{cor:MPR-cont-W1p-embed} and~\ref{cor:MPR-cont-W1q-embed}---and thus formulated in sufficient generality.

\begin{lemma}
 \label{lem:complex-interpol-in-fracpower-domain}
  Let $\zeta < 1 - \frac2s$ and $0 \leq \beta < 1-\frac2s - \zeta$. Then there is $\theta < 1-\frac1s$ such that
  \begin{equation*}
\bigl[H^{-\zeta,p}_D,\Dkal_\zeta\bigr]_\theta \hookrightarrow \domLp^{1/2+\beta/2}.
  \end{equation*}
  In particular, for $0 \leq \alpha < \frac12(1-\frac2s - \zeta-\beta)$,
  \begin{equation*}
 \W^{1,s}(I,H^{-\zeta,p}_D,\Dkal_\zeta) \hookrightarrow_c C^\alpha\bigl(\overline I,\domLp^{1/2+\beta/2}\bigr).
  \end{equation*}
\end{lemma}

\begin{proof}
  We observe that due to the assumptions on $\zeta$ and $\beta$, we can choose $\theta$ such that $\frac12(1 + \zeta + \beta) < \theta < 1-\frac1s$. Thus, with Lemma~\ref{lem:interpolation-identify-frac-power},
\begin{equation*}
   \bigl[H^{-\zeta,p}_D,\Dkal_\zeta\bigr]_{\theta} \hookrightarrow \bigl[H^{-\zeta,p}_D,\Dkal_\zeta\bigr]_{\frac12(1 + \zeta + \beta)} = \domLp^{1/2+\beta/2},
  \end{equation*}
  which is the first assertion. The second one follows immediately using Lemma~\ref{lem:max-reg-embeddings} and letting $\theta \downarrow \frac12(1 + \zeta + \beta)$.
\end{proof}

Note that Lemma~\ref{lem:complex-interpol-in-fracpower-domain} and the Kato square root property guarantee $\domLp^{1/2} = W^{1,p}_D$ for $\beta=0$. We collect the implications in the following Corollary~\ref{cor:MPR-cont-W1p-embed}.

\begin{corollary}\label{cor:MPR-cont-W1p-embed}
  Suppose that $\zeta < 1 - \frac2s$. Then for $0 \leq \alpha < \frac12(1-\frac2s - \zeta)$, we have
\begin{equation*}
   \W^{1,s}(I,H^{-\zeta,p}_D,\Dkal_\zeta) \hookrightarrow_c C^\alpha\bigl(\overline I,W^{1,p}_D\bigr).
\end{equation*}
\end{corollary}
To eventually characterize the fractional power domains $\domLp^{1/2+\beta/2}$ in
the case $\beta>0$ of Lemma~\ref{lem:complex-interpol-in-fracpower-domain} further, at least in
the higher regularity setting of Assumption~\ref{assu:regularity-settings}, we show that gradients of elements of that domain are in fact in a fractional Sobolev space. As an intuition, from the Kato square root property~\eqref{eq:kato-general}, we do know that $\nabla \Dom(A+1)^{1/2} \subseteq L^p$, so a higher order fractional power domain will provide additional regularity for gradients of its elements.

\begin{lemma}\label{lem:regular-case-fracpower-Lipschitz} Suppose that $\domLp = W^{2,p} \cap W^{1,p}_D$. Then
\begin{equation*}
\nabla \colon \domLp^{1/2+\beta/2} \to H^{\beta,p}
\end{equation*} is a continuous linear operator for every $0 \leq \beta \leq 1$.
\end{lemma}

\begin{proof}
  Note that due to the optimal regularity assumption and the Kato square root property as in Proposition~\ref{prop:kato}, we have
  \begin{equation*}
\nabla (A+1)^{-\frac12} \in \Lkal(L^p) \cap \Lkal(W^{1,p}_D,W^{1,p}).
  \end{equation*}
  Indeed, continuity of the $L^p$ version of $\nabla (A+1)^{-\frac12}$ is part of the Kato square root property. On the other hand, $(A+1)^{-\frac12}$ maps $W^{1,p}_D =\domLp^{1/2}$ continuously into $\domLp = W^{2,p} \cap W^{1,p}_D$, so $\nabla (A+1)^{-\frac12}$ is also a continuous linear operator taking $W^{1,p}_D$ to $W^{1,p}$. Interpolating between $L^p$ and $W^{1,p}$ and using Lemma~\ref{lem:interpolation-identify-frac-power}, we find that
  \begin{equation*}
 \nabla (A+1)^{-\frac12} \in \Lkal\bigl(\domLp^{\beta/2},H^{\beta,p}\bigr) \qquad \text{for all}~0 \leq \beta \leq 1,
  \end{equation*}
  and so via~\eqref{eq:fracpower-equals-bessel-below-12}, for all $y \in \domLp^{1/2+\beta/2}$,
  \begin{equation*}
\norm{\nabla y}_{H^{\beta,p}} \leq C \norm{(A+1)^{1/2}y}_{\domLp^{\beta/2}} = C \norm{y}_{\domLp^{1/2+\beta/2}}.\qedhere
  \end{equation*}
\end{proof}

We put Lemma~\ref{lem:regular-case-fracpower-Lipschitz} to good use in the next two results. In the following one, we exclude the case $p = q$ since it is already covered by Corollary~\ref{cor:MPR-cont-W1p-embed}.

\begin{corollary}\label{cor:MPR-cont-W1q-embed}
  Suppose that $\domLp = W^{2,p} \cap W^{1,p}_D$. Let $p < q \leq \infty$ and $\zeta < 1 - \frac2s - \frac{d}p + \frac{d}q$. Then for $0 \leq \alpha < \frac12(1-\frac2s - \frac{d}p + \frac{d}q - \zeta)$, we have
\begin{equation*}
   \W^{1,s}(I,H^{-\zeta,p}_D,\Dkal_\zeta) \hookrightarrow_c C^\alpha\bigl(\overline I,W^{1,q}_D\bigr)
\end{equation*}
if $q< \infty$, and for $q = \infty$:
\begin{equation*}
   \W^{1,s}(I,H^{-\zeta,p}_D,\Dkal_\zeta) \hookrightarrow_c C^\alpha\bigl(\overline I,C^1(\overline\Omega)\bigr).
\end{equation*}
\end{corollary}

\begin{proof}
   We use Lemmata~\ref{lem:complex-interpol-in-fracpower-domain} and~\ref{lem:regular-case-fracpower-Lipschitz}. Pick $\beta$ between $\frac{d}p - \frac{d}q$ and $1 - \frac2s - \zeta$. Then, Lemma~\ref{lem:complex-interpol-in-fracpower-domain} shows that
  \begin{equation*}
   \W^{1,s}(I,H^{-\zeta,p}_D,\Dkal_\zeta) \hookrightarrow_c C^\alpha\bigl(\overline I,\domLp^{1/2+\beta/2}\bigr)
\end{equation*}
for $0 \leq \alpha < \frac12(1-\frac2s-\zeta - \beta)$. Approaching $\beta \downarrow \frac{d}p - \frac{d}q$, we obtain the desired open range for $\alpha$. It remains to show that for any $\beta > \frac{d}p - \frac{d}q$, we have $\domLp^{1/2+\beta/2} \hookrightarrow W^{1,q}_D$ if $q < \infty$ and $\domLp^{1/2+\beta/2} \hookrightarrow C^1(\overline\Omega)$ if $q = \infty$. But this follows precisely by the condition that $\beta > \frac{d}p - \frac{d}q$, since it implies that $H^{\beta,p} \hookrightarrow L^q$ if $q < \infty$, and $H^{\beta,p} \hookrightarrow C(\overline\Omega)$ if $q = \infty$. The claim follows via Lemma~\ref{lem:regular-case-fracpower-Lipschitz} and $W^{1,p}_D \cap W^{1,q} = W^{1,q}_D$ for $q<\infty$.
\end{proof}

In Corollary~\ref{cor:MPR-cont-W1q-embed}, to be very precise, for $q = \infty$ we would even obtain $C^\alpha(\overline I,C^1_D(\overline\Omega))$ in the sense of strong zero boundary values on $D$, and also some degree of H\"older continuity on the gradients. However, we will not need this particular information explicitly later on.

The next lemma is the essential piece that facilitates all the maximal parabolic regularity arguments we are going to employ in Lemma~\ref{lem::mprhminuszeta} below. Note that there is, in contrast to Lemma~\ref{lem:complex-interpol-in-fracpower-domain} and its corollaries, no direct assumption on $\zeta$ here, only in relation to the degree of continuity.

\begin{lemma}\label{lem:maxregembed-hoelder-in-space}
There are $\theta < 1-\frac1s$ and $\sigma > 1-\zeta$, with $\sigma = 1$ if $\zeta = 0$, such that
\begin{equation}
  \bigl[H^{-\zeta,p}_D,\Dkal_\zeta\bigr]_\theta \hookrightarrow C^\sigma(\overline\Omega).\label{eq:interpol-space-hoelder}
\end{equation}
In particular,
\begin{equation*}
   \W^{1,s}(I,H^{-\zeta,p}_D,\Dkal_\zeta) \hookrightarrow_c C\bigl(\overline I,C^\sigma(\overline\Omega)\bigr).
\end{equation*}
\end{lemma}

\begin{proof} Before we start with the detailed proof, let us note that it will be enough to show the first assertion; due to $\theta < 1-\frac1s$, the second assertion follows immediately using~\eqref{eq:maxreg-embed}. Moreover, throughout the proof, we will consider several choices of $\theta$, all of which will satisfy $\frac12 + \frac{d}{2p} < \theta < 1-\frac1s$. Then $\theta > \frac12 > \frac\zeta2$, the latter allowing to use~\eqref{eq:interpolation-identify-frac-power}, so that it is sufficient to show the first assertion for $\domLp^{\theta-\zeta/2}$ instead of the interpolation space.

Now, for the actual proof, we branch along the different cases of Assumption~\ref{assu:regularity-settings} and then along the possible relations between $\theta$ and $\zeta$. Suppose first that $\zeta > \frac{d}p$. Then we can pick
 $\theta$ such that $\frac12 + \frac{d}{2p} < \theta < (\frac12 + \frac\zeta2) \wedge (1-\frac1s)$ which implies that $\sigma \coloneqq 2\theta - \zeta-\frac{d}p > 1-\zeta$ and $\theta - \frac\zeta2 < \frac12$. Thus, using Lemma~\ref{lem:regular-case-interpolation} and Sobolev embedding, we find that, as desired,
 \begin{equation*}
   \domLp^{\theta-\zeta/2} = H^{2\theta-\zeta,p}_D \hookrightarrow C^\sigma(\overline\Omega).
 \end{equation*}

If $\domLp = {W^{2,p} \cap W^{1,p}_D}$, then we can dispose of the requirement that $\zeta > \frac{d}p$. Indeed, let $\frac12 + \frac{d}{2p} < \theta < 1-\frac1s$. If $\theta-\frac\zeta2 \leq \frac12$, then repeat the proof as before. If not, then proceed as follows. Define $\gamma > 0$ via $\theta-\frac\zeta2 = \frac12 + \frac\gamma2$.

First suppose that $\theta - \frac\zeta2 > \frac12 + \frac{d}{2p}$. Then $\gamma > \frac{d}{p}$ and Lemma~\ref{lem:regular-case-fracpower-Lipschitz} shows that $\nabla z \in H^{\gamma,p} \hookrightarrow C(\overline\Omega)$ for every $z \in \domLp^{1/2 + \gamma/2}$. It follows that
\begin{equation*}
  \domLp^{\theta - \zeta/2} = \domLp^{1/2 + \gamma/2} \hookrightarrow C^{1}(\overline\Omega),
\end{equation*}
which completes the proof for this case, since clearly $\sigma = 1 > 1-\zeta$, including the special variant $\sigma = 1$ for $\zeta = 0$.

Finally, suppose that $\frac12 < \theta - \frac\zeta2 \leq \frac12 + \frac{d}{2p}$. Then $0 < \gamma \leq \frac{d}{p}$. Again using Lemma~\ref{lem:regular-case-fracpower-Lipschitz}, we obtain that $\nabla z \in H^{\gamma,p} \hookrightarrow L^q$ for $q$ defined by $\frac{d}q = \frac{d}p - \gamma$ if $\gamma < \frac{d}p$, and for every $q < \infty$ for $\gamma = \frac{d}p$, for every $z \in \domLp^{1/2 + \gamma/2}$. It follows that
\begin{equation*}
  \domLp^{\theta - \zeta/2} = \domLp^{1/2 + \gamma/2} \hookrightarrow W^{1,q} \hookrightarrow C^\sigma(\overline\Omega),
\end{equation*}
with $\sigma = 1-\frac{d}q = 1 + \gamma - \frac{d}p = 2\theta - \zeta - \frac{d}p$, if $\gamma < \frac{d}p$, and that indeed satisfies $\sigma > 1-\zeta$ due to the initial choice of $\theta > \frac12 + \frac{d}{2p}$. If $\gamma = \frac{d}p$, then we can choose $q$ as large as we desire, in particular, we can choose it so that $\frac{d}q < \zeta$ and so $\sigma = 1-\frac{d}q > 1-\zeta$. This completes the proof.
\end{proof}

\begin{remark}
 In the regular case $\domLp = W^{2,p} \cap W^{1,p}_D$, the second assertion in Lemma~\ref{lem:maxregembed-hoelder-in-space} can also be obtained directly from Corollary~\ref{cor:MPR-cont-W1q-embed}: Choose $q>p$ such that $0 < \zeta - \frac{d}q < 1-\frac2s - \frac{d}p$ and observe that $W^{1,q} \hookrightarrow C^{\sigma}(\overline\Omega)$ with $\sigma\coloneqq 1 - \frac{d}q > 1-\zeta$ by construction. However, we will require~\eqref{eq:interpol-space-hoelder} explicitly later on, so we cannot dispose of the slightly fiddly proof of Lemma~\ref{lem:maxregembed-hoelder-in-space}.
\end{remark}

During the next chapter we will repeatedly make use of the following result on superposition operators that follows from a short straight forward computation:

\begin{lemma}\label{lem::superposition}
For any $0 \leq \sigma \leq 1$, the function $\xi$ induces a continuously Fréchet-differentiable superposition operator $\Xi \colon C^\sigma(\overline\Omega) \rightarrow C^\sigma(\overline\Omega)$ with the derivative $\Xi'(y)z = \xi'(y)z$ for $y,z \in C^\sigma(\overline\Omega)$.
\end{lemma}

The following results on maximal parabolic regularity will be vital for the analysis of the state equation, both for local-in-time wellposedness and the optimal control context in Sections~\ref{sec:global-control-prop},~\ref{sec::existence} and~\ref{sec::fon}. For the definitions of $\Akal$ and $\Alin$, see~\eqref{eq:lin-A-def}; see also Proposition~\ref{prop:A-mpr-first-statement}.

\begin{lemma}\label{lem::mprhminuszeta}
Fix $y \in \W^{1,s}(I,H^{-\zeta,p}_D,\Dkal_\zeta)$.
Then we have the following:
\begin{enumerate}[(i)]
\item For every $r \in (1,\infty)$, the operator $\Akal(y)$ satisfies nonautonomous maximal parabolic regularity on $L^r(I,H^{-\zeta,p}_D)$ with constant domain $\Dkal_\zeta$, and $t \mapsto \Akal(y)(t)$ is continuous as a mapping $I \rightarrow \Lkal(\Dkal_\zeta,H^{-\zeta,p}_D)$.
\item The operator $\Akal(y) + \Alin(y)$ satisfies nonautonomous maximal parabolic regularity on $L^s(I,H^{-\zeta,p}_D)$ with constant domain $\Dkal_\zeta$.
\item Let Assumption~\ref{gassemilin}~\ref{gassemilin:diffbar} hold true. Then the operator $\Akal(y) + \Alin(y) - \Fkal'(y)$ satisfies nonautonomous maximal parabolic regularity on $L^s(I,H^{-\zeta,p}_D)$ with constant domain $\Dkal_\zeta$.
\end{enumerate}
All assertions stay true when $I$ is replaced by any subinterval of $I$ \emph{mutatis mutandis}.
\end{lemma}

In the following proof, for $y \in \W^{1,s}(I,H^{-\zeta,p}_D,\Dkal_\zeta)$ and $v$ not depending on $t$, we will freely use the identification $[\Alin(y)v](t) \simeq [\Alin(y)(t)]v$ and analogously for $\Fkal'(y)$.

\begin{proof}[Proof of Lemma~\ref{lem::mprhminuszeta}]
   Let $y \in \W^{1,s}(I,H^{-\zeta,p}_D,\Dkal_\zeta)$. We rely on~\cite[Lemma~4.2/Theorem~4.3]{HoppeMeinlschmidtNeitzel2023} which gives conditions for nonautonomous maximal parabolic regularity of $\Akal$ in $H^{-\zeta,p}_D$ with constant domain $\Dkal_\zeta$ building on H\"older continuity in space and continuity in time of the perturbative coefficient $\xi(y)$. By this technique, nonautonomous maximal parabolic regularity on every subinterval of $I$ comes for free. To this end, recall that $\W^{1,s}(I,H^{-\zeta,p}_D,\Dkal_\zeta) \hookrightarrow C(\overline I,C^\sigma(\overline\Omega))$ for some $\sigma > 1-\zeta$, and $\sigma =1$ if $\zeta = 0$, as established in Lemma~\ref{lem:maxregembed-hoelder-in-space}. We note that in~\cite{HoppeMeinlschmidtNeitzel2023}, we consider $\Akal(y)+1$ instead of $\Akal(y)$---which is why $\Dkal_\zeta$ is also defined using ``$+1$''---but it is easily seen that the maximal parabolic regularity results for the former transfer immediately to the latter with an $e^{-t}$ rescaling of the associated solutions.

   \begin{enumerate}[(i)]
   \item For $\Akal(y)$, we first use~\cite[Lemma~4.2]{HoppeMeinlschmidtNeitzel2023} which says that since $y\in C(\overline I,C^\sigma(\overline\Omega))$, the function $t \mapsto \Akal(y)(t)$ is continuous $I \rightarrow \Lkal(\Dkal_\zeta,H^{-\zeta,p}_D)$; as a consequence~(\cite[Theorem~4.3]{HoppeMeinlschmidtNeitzel2023}---see also Proposition~\ref{prop:cont-mpr-plus-pruess})---$\Akal(y)$ satisfies nonautonomous maximal parabolic regularity on $L^r(I,H^{-\zeta,p}_D)$ with constant domain $\Dkal_\zeta$ for every $r \in (1,\infty)$.

  \item Nonautonomous maximal parabolic regularity of $\Akal(y) + \Alin(y)$ on $L^s(I,H^{-\zeta,p}_D)$ is then obtained from the perturbation result of Pr\"uss as in Proposition~\ref{prop:cont-mpr-plus-pruess}. In the present case, verifying the assumptions of that perturbation result amounts to showing that $\Alin(z) \in L^s(I,\Lkal((H^{-\zeta,p}_D,\Dkal_\zeta)_{1-1/s,s},H^{-\zeta,p}_D))$. We will in fact prove that $\Alin(y) \in L^s(I,\Lkal(C^\sigma(\overline\Omega),H^{-\zeta,p}_D))$ whenever $\sigma > 1-\zeta$. Due to the standard embedding $(H^{-\zeta,p}_D,\Dkal_\zeta)_{1-1/s,s} \hookrightarrow [H^{-\zeta,p}_D,\Dkal_\zeta]_\theta$ for any $\theta < 1-\frac1s$ and Lemma~\ref{lem:maxregembed-hoelder-in-space}, this will be sufficient.

  In fact, let $\sigma>1-\zeta$ be the number from Lemma~\ref{lem:maxregembed-hoelder-in-space}, and let $v \in C^\sigma(\overline\Omega)$. With the assumptions on $\xi$ as in~\eqref{gaseq}, we then have $\xi'(y(t))v \in C^\sigma(\overline\Omega)$ for every $t \in I$, and $t \mapsto \xi'(y(t))v \in C(\overline I,C^\sigma(\overline\Omega))$. Again relying on~\cite[Lemma~4.2]{HoppeMeinlschmidtNeitzel2023}, this implies that $-\nabla \cdot \xi'(y)v \mu\nabla \in C(\overline I,\Lkal(\Dkal_\zeta,H^{-\zeta,p}_D))$ with
  \begin{equation*}
\sup_{t\in I} \norm{-\nabla \cdot \xi'(y(t))v \mu\nabla}_{\Lkal(\Dkal_\zeta,H^{-\zeta,p}_D)} \leq \sup_{t \in I} \norm{\xi'(y(t))}_{C^\sigma(\overline\Omega)} \, \norm{v}_{C^\sigma(\overline\Omega)}.
  \end{equation*}
  It follows that for $y \in \W^{1,s}(I,H^{-\zeta,p}_D,\Dkal_\zeta)$ and $t \in I$,
  \begin{align*}
\norm[\big]{[\Alin(y)(t)]v}_{H^{-\zeta,p}_D} & = \norm{-\nabla \cdot \xi'(y(t))v \mu\nabla y(t)}_{H^{-\zeta,p}_D} \\ & \leq \norm{y(t)}_{\Dkal_\zeta} \, \norm{v}_{C^\sigma(\overline\Omega)} \, \sup_{t \in I} \norm{\xi'(y(t))}_{C^\sigma(\overline\Omega)}.
  \end{align*}
  We obtain
  \begin{equation*}
\norm[\big]{\Alin(y)(t)}_{\Lkal(C^\sigma(\overline\Omega),H^{-\zeta,p}_D)} \leq \norm{y(t)}_{\Dkal_\zeta} \, \sup_{t \in I} \norm{\xi'(y(t))}_{C^\sigma(\overline\Omega)},
  \end{equation*}
  whose right-hand side as a function of $t$ is in $L^s(I)$ since $y \in L^s(I,\Dkal_\zeta) \cap C(\overline I,C^\sigma(\overline\Omega))$. That concludes the proof.

  \item The third assertion follows immediately from the foregoing proof, since by Assumption~\ref{gassemilin}~\ref{gassemilin:diffbar} and what was proven before,
  \begin{equation*}
\Alin(y) - \Fkal'(y) \in L^s\bigl(I,\Lkal((H^{-\zeta,p}_D,\Dkal_\zeta)_{1-1/s,s},H^{-\zeta,p}_D)\bigr)
  \end{equation*}
  is an admissible perturbation for continuous-in-time nonautonomous maximal parabolic regularity as offered by $\Akal(z)$; see Proposition~\ref{prop:cont-mpr-plus-pruess}.\qedhere
  \end{enumerate}
\end{proof}

\section{Analysis of solutions of the state equation}\label{sec::solutions}

In this section we discuss the state equation~\eqref{se} as an abstract parabolic evolution equation in $H^{-\zeta,p}_D$. For this, we, again, presume that Assumptions~\ref{ass:geometry},~\ref{ass:coeff-matrix},~\ref{ass:optimal-elliptic-regularity-mu},~\ref{ass:xi},~\ref{assu:regularity-settings} and~\ref{gaseq} are satisfied without further ado, and we will only state which parts of Assumption~\ref{gassemilin} we require for the respective results. The basic relations on integrabilites from these assumptions are $p>d$ and $\frac2s < 1-\frac{d}p$. There are no particular restrictions on the scale parameter $\zeta$, however, we point out that in the case of irregular data as in Assumption~\ref{assu:regularity-settings}, we must stick with $\zeta > \frac{d}p$. We begin with the proof of existence of local-in-time solutions, invoking a classical result of Amann; cf.~\cite[Theorem~2.1]{Amann2005_3} or~\cite{Amann2005_4}.
Since we do not pose any kind of \emph{global} growth condition on the nonlinearity $\Fkal$ (or a sign condition), we cannot expect global-in-time solutions to~\eqref{se} at all at this stage. Yet, of course, in the context of the optimal control problem~\eqref{ocp}, we are certainly interested in global-in-time solutions to~\eqref{se} which we intend to achieve by the interplay of the gradient constraints and the \emph{particular} growth conditions on $\Fkal$ as in Assumption~\ref{gassemilin}~\ref{gassemilin:boundedness}. We thus pave the way for these arguments by discussing sufficient conditions for global-in-time existence of given solutions in Section~\ref{sec:global-in-time}. A starting point here are our recent results in~\cite{HoppeMeinlschmidtNeitzel2023} which guarantee global-in-time existence of solutions if $\Fkal$ does in fact not depend on the searched-for function at all; see Theorem~\ref{thm::existenceBoNe} below. We conclude this section with a discussion of the set of \emph{global controls}, that is, controls $u$ for which the associated solution to~\eqref{se} in fact exists globally in time, and it will turn out that this is an \emph{open} set under suitable differentiability assumptions on the nonlinear functions involved in~\eqref{se}, but it is also weakly closed in particular situations. Such properties will prove most useful when establishing existence of globally optimal solutions to~\eqref{ocp} in Section~\ref{sec::existence} and also for the derivation of first order necessary optimality conditions for~\eqref{ocp} in Section~\ref{sec::fon}.

Let us point out that the approach to enforce global-in-time solutions to the state equation by leveraging constraints to the optimal control problem borrows several ideas from~\cite{Meinlschmidt2017_1,Meinlschmidt2017_2}.

\subsection{Existence and uniqueness of local-in-time solutions}\label{sec::local-in-time}

First, we recall a result on wellposedness and global-in-time existence of solutions to~\eqref{se} when $\Fkal$ does in fact not depend on $y$; see~\cite[Theorem~5.7, also Remark~5.9]{HoppeMeinlschmidtNeitzel2023}. It is a generalization of the results in~\cite{Bonifacius2018} which are in turn based on~\cite{Meinlschmidt2016}. This particular theorem will be highly convenient for the following, as it essentially decouples possible blow up phenomena induced by the quasilinear nonlinearity---which do not exist in the present case, as the theorem states---from the ones induced by $\Fkal$.

\begin{theorem}\label{thm::existenceBoNe} For any $f \in L^s(I,H^{-\zeta,p}_D)$ and $y_0 \in (H^{-\zeta,p}_D,\Dkal_\zeta)_{1-\frac1s,s}$, there exists a unique solution $y \in \W^{1,s}(I,H^{-\zeta,p}_D,\Dkal_\zeta)$ to the quasilinear parabolic equation
  \begin{equation*}
\begin{aligned}
\partial_t y
+ \Akal(y)y &= f && \text{ in } H^{-\zeta,p}_D \text{ a.e.\ on } I, \\
y(0) &= y_0 && \text{ in } (H^{-\zeta,p}_D,\Dkal_\zeta)_{1-\frac1s,s}.
\end{aligned}
\end{equation*}
Moreover, the associated solution operator $(f,y_0) \mapsto z$ maps bounded sets into bounded sets.
\end{theorem}

In the present case of the state equation~\eqref{se} with $\Fkal$ not admitting any global growth bounds, however, we cannot reasonably expect that a general global-in-time existence result can be established; cf.\ the discussion of the (counter) examples~\eqref{eq:counterexF1} and~\eqref{eq:counterexF2} in the introduction.
We will therefore proceed differently in the following. Indeed, we accept that we can only achieve local-in-time solutions to~\eqref{se} in the present context, however, we leverage the bounds induced by the gradient constraints to obtain a well posed optimal control problem~\eqref{ocp}.

Let us thus recall the concept of maximal local-in-time and global-in-time solutions for~\eqref{se}.
\begin{definition}[Maximal and global solution]\label{def:maximal-global-solution}
 Given $u \in L^s(\Lambda)$, we say that $y$ a \emph{maximal} solution to~\eqref{se} if there is $T^\bullet_u \in (0,T]$ such that for every interval $J = (0,\tau)$ with $\tau \in (0,T^\bullet_u)$, the following are true:
\begin{itemize}
\item $y \in \W^{1,s}(J,H^{-\zeta,p}_D,\Dkal_\zeta)$ with $y(0)=y_0$ in $(H^{-\zeta,p}_D,\Dkal_\zeta)_{1-\frac1s,s}$,
\item $y$ is a solution to~\eqref{se} on $J$, that is, for almost all $t \in J$,
\[
   \partial_t y(t) + \Akal(y(t))y(t) = \Fkal(y)(t) + (Bu)(t) \quad \text{in}~H^{-\zeta,p}_D,
\]
\item the number $T^\bullet_u \in (0,T]$ is the \emph{maximal} one having the foregoing properties.
\end{itemize}
  If the foregoing assertions in fact hold true for $J = [0,T]$, then we say that $y$ is a \emph{global} solution to~\eqref{se}. Otherwise, $T^\bullet_u$ is the \emph{maximum time of existence} for $y$.
\end{definition}

The following result guarantees unique maximal local-in-time solutions for~\eqref{se} and will serve as a basis for all further considerations.

\begin{proposition}\label{prop::localex}
Let Assumption~\ref{gassemilin}~\ref{gassemilin:volterra} on $\Fkal$ and the remaining assumptions of Section~\ref{sec:notass} hold. Then, for any $u \in L^s(\Lambda)$, there
is a unique maximal local-in-time solution to~\eqref{se} in the sense of Definition~\ref{def:maximal-global-solution}.
\end{proposition}

\begin{proof}
We use the seminal theorem by Amann, cf.~\cite[Theorem~2.1]{Amann2005_3} and also~\cite{Amann2005_4}, with $E_0 = H^{-\zeta,p}_D$ and $E_1 = \Dkal_\zeta$. The required properties of $\Fkal$ are covered precisely by Assumption~\ref{gassemilin}~\ref{gassemilin:volterra}, and the initial value $y_0$ has the correct regularity due to Assumption~\ref{assu:regularity-settings}. Regarding the quasilinear operator $\Akal$, we can rely on Lemma~\ref{lem::mprhminuszeta}, which gives all the required properties. The Volterra property for $\Akal$ is evident, so the cited theorem states precisely the assertion.
\end{proof}

\subsection{Global-in-time-solutions}\label{sec:global-in-time} In this section we consider the intricate question in which situations the maximal-in-time solutions to~\eqref{se} provided by Proposition~\ref{prop::localex} are in fact global-in-time. In view of the optimal control problem~\eqref{ocp} under consideration, this is most important to understand, since both the cost function in~\eqref{ocp} as well as the constraints on the (gradient of the) state $y$ as in~\eqref{e:yadavergd} and~\eqref{e:yadptwise} explicitly require $y$ to be defined on the \emph{whole} given time interval.

First, we prove that maximal local-in-time solutions where the right-hand sides exist ``globally'' are in fact global-in-time solutions in the sense of Definition~\ref{def:maximal-global-solution}, so, we derive a sufficient condition. Further, and more importantly for what follows, we derive both openness and a suitable weak closedness property for the set of controls that allow for global-in-time solutions, in essence following arguments from~\cite{Meinlschmidt2017_1}.

\begin{proposition}\label{prop::globalex}
Let Assumption~\ref{gassemilin}~\ref{gassemilin:volterra} on $\Fkal$
hold true, and let $y$ be the maximal-local-in-time solution of~\eqref{se} for given $u \in L^s(\Lambda)$ with maximal time of existence $T_u^\bullet$ in the sense of Definition~\ref{def:maximal-global-solution}.
Suppose that in fact $\Fkal(y) \in L^s(0,T^\bullet_u,H^{-\zeta,p}_D)$.
Then $T_u^\bullet =T$ and $y \in \W^{1,s}(I,H^{-\zeta,p}_D,\Dkal_\zeta)$, that is, $y$ is the global-in-time solution to~\eqref{se}. Further, there is a constant $C > 0$ depending only on the problem data, the norm of $\Fkal(y)$ in $L^s(I,H^{-\zeta,p}_D)$ and $\lVert u \rVert_{L^s(\Lambda)}$ such that $\lVert y \rVert_{\W^{1,s}(I,H^{-\zeta,p}_D, \Dkal_\zeta)} \leq C$.
\end{proposition}
Before we start the proof, let us point out that a growth bound on $\Fkal$ as in Assumption~\ref{gassemilin}~\ref{gassemilin:boundedness} is not required in the foregoing proposition. Indeed it is an \emph{assumption} that $\Fkal(y) \in L^s(0,T^\bullet_u,H^{-\zeta,p}_D)$. 
Eventually, in Section~\ref{sec::existence}, Proposition~\ref{prop::globalex} is indeed leveraged by~Assumption~\ref{gassemilin}~\ref{gassemilin:boundedness} via Proposition~\ref{prop::closedness} to obtain weak closedness for the admissible sets in the control problems. 

\begin{proof}[Proof of Proposition~\ref{prop::globalex}]
Observe that by assumption, $Bu + \Fkal(y) \in L^s(0,T^\bullet_u,H^{-\zeta,p}_D)$ with
\begin{multline}\label{eq0}
\lVert Bu + \Fkal(y) \rVert_{L^s(0,T^\bullet_u,H^{-\zeta,p}_D)} \leq \lVert B \rVert_{\Lkal(L^s(\Lambda), L^s(I,H^{-\zeta,p}_D))} \lVert u \rVert_{L^s(\Lambda)} \\+\lVert\Fkal(y)\rVert_{L^s(0,T^\bullet_u,H^{-\zeta,p}_D)}.
\end{multline}
Therefore, by Theorem~\ref{thm::existenceBoNe} there is a unique solution $w \in \W^{1,s}(0,T^\bullet_u,H^{-\zeta,p}_D,\Dkal_\zeta)$ to the equation
\begin{align*}
\begin{aligned}
\partial_t w
+ \Akal(w)w &= Bu + \Fkal(y)& &\text{in } H^{-\zeta,p}_D \text{ a.e.\ on } (0,T^\bullet_u) \\
w(0) &= y_0 &&\text{in } (H^{-\zeta,p}_D,\Dkal_\zeta)_{1-\frac1s,s}.
\end{aligned}
\end{align*}
With uniqueness of local-in-time solutions, cf.\ Proposition~\ref{prop::localex}, we find that $w \equiv y$ on $[0,T^\bullet_u)$. In particular, $y \in \W^{1,s}(0,T^\bullet_u,H^{-\zeta,p}_D,\Dkal_\zeta)$, and this already implies that $T^\bullet_u = T$ and that $y$ is the global-in-time solution to~\eqref{se}; cf.~\cite[Corollary~3.2]{Pruess2002} or~\cite[Theorem~2.1]{Amann2005_3}.

Further, since the right-hand side in~\eqref{eq0} only depends on the problem data, $\lVert u \rVert_{L^s(\Lambda)}$ and $\lVert\Fkal(y)\rVert_{L^s(0,T^\bullet_u,H^{-\zeta,p}_D)}$, the claimed bound for $ \lVert y \rVert_{\W^{1,s}(I,H^{-\zeta,p}_D,\Dkal_\zeta)}$ also follows from Theorem~\ref{thm::existenceBoNe}.
\end{proof}

In view of the optimal control problem~\eqref{ocp}, we summarize that Proposition~\ref{prop::localex} provides existence of local-in-time-solutions to the state equation for all admissible controls $u\in U_{\ad}$, and Proposition~\ref{prop::globalex} postulates a sufficient condition that leads to global solutions. We now introduce the set of \emph{global controls} $U_g \subseteq L^s(\Lambda)$ which admit \emph{global-in-time} solutions to~\eqref{se}. On this set, we have a well defined control-to-state operator
\begin{equation}\label{def:S}
S\colon U_g \rightarrow \W^{1,s}(I,H^{-\zeta,p}_D, \Dkal_\zeta), \qquad u \mapsto y=y(u),
\end{equation}
that maps a control $u \in U_g$ to the state $y=S(u)$ associated to $u$ via~\eqref{se}. The fact that $y = S(u)$ is a \emph{global-in-time} solution to~\eqref{se} is encoded in the function space associated with the whole time interval $I$ in~\eqref{def:S}. We now establish several properties of $U_g$ and $S$.

\subsection{Properties of the set \texorpdfstring{$U_g$}{of global controls}}\label{sec:global-control-prop}
As mentioned in the introduction of this section, we will now show auxiliary properties for the proof of existence of optimal controls and the derivation of first order necessary optimality conditions.
We first prove that $U_g$ is open and that the control-to-state map $S$ is differentiable.

\begin{proposition}\label{prop::openset}
Let Assumptions~\ref{gassemilin}~\ref{gassemilin:volterra} and~\ref{gassemilin:diffbar} be satisfied. Then the following assertions hold true:
\begin{enumerate}[(i)]
\item The set $U_g \subset L^s(\Lambda)$ is open.
\item The control-to-state map $S\colon U_g \rightarrow \W^{1,s}(I,H^{-\zeta,p}_D, \Dkal_\zeta)$ is continuously Fréchet differentiable and its derivative $S'(u) v = z$ in direction $v$ is given by the unique solution $z$ to the \emph{sensitivity equation}
\begin{align}\label{lse}\tag{SEq}
\begin{aligned}
\partial_t z + \Akal(y) z + \Alin(y)z - \Fkal'(y) z &= B v && \text{in } H^{-\zeta,p}_D \text{ a.e.\ on } I, \\
z(0) &= 0 && \text{in } (H^{-\zeta,p}_D, \Dkal_\zeta)_{1-\frac1s,s}.
\end{aligned}
\end{align}
\end{enumerate}
\end{proposition}

\begin{proof}
The proof is similar to~\cite[Theorem 3.1]{Meinlschmidt2017_2}.
If $U_g = \emptyset$, then there is nothing to show. Otherwise, fix an arbitrary $\hat u \in U_g$ and its associated state $\hat y\coloneqq S(\hat u)$. We apply the implicit function theorem to
\begin{align*}
E\colon \W^{1,s}(I,H^{-\zeta,p}_D,\Dkal_\zeta) \times L^s(\Lambda) &\rightarrow L^s(I,H^{-\zeta,p}_D) \times (H^{-\zeta,p}_D,\Dkal_\zeta)_{1-1/s,s}, \\
(y,u) &\mapsto \bigl( \partial_t y + \Akal(y)y - Bu - \Fkal(y), y(0) - y_0 \bigr)
\end{align*}
at $(y,u) = (\hat y, \hat u)$. This will imply both assertions since the control-to-state operator $u \mapsto S(u)$ is exactly the implicit function realizing $E(S(u),u) = 0$ for all $u \in U_g$. Indeed, by construction, we have $E(\hat y, \hat u) = 0$. The given assumptions on $\xi$ and $\Fkal$ imply that the nonlinear operator $E$ is continuously differentiable. This follows from~\cite[Theorems~4 and~7]{Goldberg1992} utilizing Lemmas~\ref{lem:maxregembed-hoelder-in-space} and~\ref{lem::superposition}, and we find
\begin{align*}
\partial_y E(\hat y,\hat u) z &= \bigl( \partial_t z + \Akal(\hat y) z + \Alin(\hat y) z - \Fkal'(\hat y) z, z(0) \bigr)\\
\partial_u E(\hat y, \hat u) v &= ( -Bv, 0),
\end{align*}
with $\Alin$ as introduced in~\eqref{eq:lin-A-def}. The partial derivative $\partial_y E(\hat y, \hat u)$ is continuously invertible as an operator $ \W^{1,s}_0(I,H^{-\zeta,p}_D,\Dkal_\zeta) \to L^s(I,H^{-\zeta,p}_D)$ precisely when the linearized differential operator $\Akal(\hat y) + \Alin(\hat y) - \Fkal'(\hat y)$ admits nonautonomous maximal parabolic regularity on $L^s(I,H^{-\zeta,p}_D)$; cf.\ Appendix~\ref{sec:mpr}. But we have already established this property in Lemma~\ref{lem::mprhminuszeta} using Assumption~\ref{gassemilin}~\ref{gassemilin:diffbar} for $\Fkal'$, so the implicit function theorem can do its magic and yield the assertion.
\end{proof}

\begin{remark}\label{rem:general-linearized-solution-operator}
Let $y \in \W^{1,s}(I,H^{-\zeta,p}_D, \Dkal_\zeta)$. In the proof of Proposition~\ref{prop::openset}, we have used that the nonautonomous parabolic operator $\Akal(y) + \Alin(y) - \Fkal'(y)$ satisfies maximal parabolic regularity on $L^s(I,H^{-\zeta,p}_D)$ with constant domain $\Dkal_\zeta$ as provided by Lemma~\ref{lem::mprhminuszeta}. We give a name to the associated inverse parabolic operator for later reference:
\begin{equation*}
  \altS(y) \coloneqq \bigl(\partial_t + \Akal(y) + \Alin(y) - \Fkal'(y)\bigr)^{-1} \colon L^s(I,H^{-\zeta,p}_D)
  \to \W^{1,s}_0(I,H^{-\zeta,p}_D, \Dkal_\zeta).
\end{equation*}
We will later relate $\altS$ with the control-to-state operator $S$ via $S'(u) = \altS(y)B$ with $y = S(u)$. In particular, $\altS(y)$ is the solution operator for the linearized state equation as in~\eqref{lse} for a general right-hand side $f$ instead of $Bv$.
\end{remark}

Proposition~\ref{prop::openset} is complemented by the following weak closedness type result for $U_g$:

\begin{proposition}\label{prop::closedness}
Let Assumptions~\ref{gassemilin}~\ref{gassemilin:volterra} and~\ref{gassemilin:weak-cont} hold, and let $(u_k) \subset U_g$ be weakly convergent in $L^s(\Lambda)$ with limit $\bar u$, such that $(\Fkal(y_k))$ is bounded in $L^s(I,H^{-\zeta,p}_D)$ for $y_k\coloneqq S(u_k)$.
Then $\bar u \in U_g$ and $y_k \rightharpoonup S(\bar u)$ in $\W^{1,s}(I,H^{-\zeta,p}_D, \Dkal_\zeta)$.
\end{proposition}

\begin{proof}
We do not relabel subsequences in this proof. Boundedness of the sequences $(u_k)$ and $(\Fkal(y_k))$
allows to apply Proposition~\ref{prop::globalex} from which we find that $(y_k)$ is bounded in $\W^{1,s}(I,H^{-\zeta,p}_D, \Dkal_\zeta)$. In particular, there is a weakly convergent subsequence of $(y_k)$ whose weak limit is $\bar y \in \W^{1,s}(I,H^{-\zeta,p}_D, \Dkal_\zeta)$. Due to weak continuity of $\Fkal$ as in Assumption~\ref{gassemilin}~\ref{gassemilin:weak-cont}, this already implies that $\Fkal(y_k) \rightharpoonup \Fkal(\bar y)$ in $L^s(I,H^{-\zeta,p}_D)$. Further, using Lemma~\ref{lem:maxregembed-hoelder-in-space}, we infer that there is a $\sigma > 1-\zeta$---with $\sigma = 1$ if $\zeta = 0$---and a subsequence of $(y_k)$ that converges to some $\bar y$ in $C(\overline I,C^\sigma(\overline\Omega))$. With~\cite[Lemma~4.2]{HoppeMeinlschmidtNeitzel2023}, it follows that $\partial_t + \Akal(y_k) \to \partial_t + \Akal(\bar y)$ in the space of bounded linear operators $\W^{1,s}(I,H^{-\zeta,p}_D, \Dkal_\zeta) \to L^s(I,H^{-\zeta,p}_D)$. But since each $\Akal(y_k)$ and $\Akal(\bar y)$ in fact satisfy nonautonomous maximal parabolic regularity on $L^s(I,H^{-\zeta,p}_D)$ with constant domain $\Dkal_\zeta$ by Lemma~\ref{lem::mprhminuszeta}, using Proposition~\ref{prop:cont-mpr-plus-pruess}, we find that
\begin{equation*}
  \bigl( \partial_t + \Akal(y_k), \delta_0 \bigr)^{-1} \to \bigl( \partial_t + \Akal(\bar y), \delta_0 \bigr)^{-1}
\end{equation*}\
in the space of linear bounded operators
\begin{equation*}
  L^s(I,H^{-\zeta,p}_D) \times (H^{-\zeta,p}_D,\Dkal_\zeta)_{1-\frac1s,s} \to \W^{1,s}(I,H^{-\zeta,p}_D, \Dkal_\zeta).
\end{equation*}
Recall from above that $y_k \to \bar y$ in $C(\overline I,C^\sigma(\overline\Omega))$, so that we obtain
\begin{equation*}
   y_k = \bigl( \partial_t + \Akal(y_k), \delta_0 \bigr)^{-1}\bigl(Bu_k + \Fkal(y_k),y_0\bigr) \rightharpoonup \bigl( \partial_t + \Akal(\bar y), \delta_0 \bigr)^{-1}\bigl(B\bar u + \Fkal(\bar y),y_0\bigr)
\end{equation*}
in $\W^{1,s}(I,H^{-\zeta,p}_D,\Dkal_\zeta)$ from the continuity of $B$ and what was established before. Thus, since the limit of $(y_k)$ must be unique in $C(\overline I,C^\sigma(\overline\Omega))$,
\begin{equation*}
  \bar y = \bigl( \partial_t + \Akal(\bar y), \delta_0 \bigr)^{-1}\bigl(B\bar u + \Fkal(\bar y),y_0\bigr),
\end{equation*}
i.e., $\bar y = S(\bar u)$ due to uniqueness of solutions to~\eqref{se}. Convergence of the full sequence $y_k \rightharpoonup S(\bar u)$ follows with a subsequence-subsequence argument.
\end{proof}

\section{The optimal control problem and existence of optimal controls}\label{sec::existence}
We now turn to the analysis of the optimal control problem~\eqref{ocp}. That is, we first recall the precise formulation of the gradient constraints introduced in~\eqref{e:yadptwise} and~\eqref{e:yadavergd}, and then we put that definition to use to state a rigorous and well defined version of~\eqref{ocp} in which we restrict the formal problem~\eqref{ocp} in an implicit manner to the set of global controls $u \in U_g$, that is, controls whose associated solution in fact exists globally in time in the sense of Definition~\ref{def:maximal-global-solution}. This takes care of the issue that both cost functional and gradient constraints require the solutions to the equation under consideration to exist globally in time.

With this understanding of the problem, we prove the first main result of this paper: Assuming existence of a feasible point, i.e., existence of an admissible control whose associated state exists globally in time and satisfies the given constraints, we obtain existence of at least one global solution to~\eqref{ocp}; cf.\ Theorem~\eqref{thm::existence}. For this section, we do not require differentiability properties of $\Fkal$ as in Assumption~\ref{gassemilin}~\ref{gassemilin:diffbar}. However, we require the setting to be a bit more regular compared to what was assumed before for wellposedness of the state equation in terms of $\zeta$, cf.~\eqref{zetaqlink} in Remark~\ref{rem:assumptions-F}, and we distinguish between different cases regarding Assumption~\ref{assu:regularity-settings}. We summarize the corresponding requirements here:
\begin{center}
  \vspace{0.75em}
  \fbox{\parbox{0.95\linewidth}{
   Suppose that Assumptions~\ref{ass:geometry},~\ref{ass:coeff-matrix},~\ref{ass:optimal-elliptic-regularity-mu},~\ref{ass:xi},~\ref{assu:regularity-settings},~\ref{gaseq},~\ref{gasconstr} and Assumption~\ref{gassemilin}~\ref{gassemilin:volterra}--\ref{gassemilin:weak-cont} hold true.
  Let $q$ be the number from Assumption~\ref{gassemilin}~\ref{gassemilin:boundedness} and suppose either one of the following settings:
  \begin{itemize}
  \item the irregular setting~\irreg{} with $\frac{d}p < \zeta < 1-\frac2s$ and $q = p$, or
  \item the regular setting~\reg{} with $\domLp = W^{2,p} \cap W^{1,p}_D$ and $\zeta$ such that $\zeta < 1-\frac2s - \frac{d}p + \frac{d}q$.
  \end{itemize}
  }}
  \vspace{0.75em}
\end{center}
Recall that $1-\frac{2}{s}>\frac{d}{p}$ due to Assumption~\ref{ass:optimal-elliptic-regularity-mu}, $p\le q$ due to Assumption~\ref{gassemilin}~\ref{gassemilin:boundedness}, and that in either case, we will have $\zeta < 1-\frac2s$ in the following, which allows to deal with our model nonlinearities~\eqref{eq:examplesF1} and~\eqref{eq:examplesF2}, including the quadratic gradient nonlinearity $\Fkal(y) = \abs{\nabla y}_2^2$; cf.~Lemma~\ref{lem:grad-nonlinearity}.
By virtue of Corollaries~\ref{cor:MPR-cont-W1p-embed} and~\ref{cor:MPR-cont-W1q-embed}, the conditions above ensure, and indeed that is their purpose, that
\begin{equation}
   \W^{1,s}(I,H^{-\zeta,p}_D,\Dkal_\zeta) \hookrightarrow_c C^\alpha\bigl(\overline I,W^{1,q}_D\bigr)\label{eq:maxreg-into-cw1q}
\end{equation}
precisely for the $q$ as in Assumption~\ref{gassemilin}~\ref{gassemilin:boundedness}, and $\Fkal$ is uniformly bounded on bounded sets in $C(\overline I,W^{1,q}_D)$. (For $q=\infty$, we replace $W^{1,q}_D$ by $C^1(\overline\Omega)$.) We will make use of this entanglement by posing gradient constraints associated precisely to the latter space, thereby facilitating a uniform boundedness of $\Fkal(y)$ over all $y$ satisfying the gradient constraints. We refer to Assumption~\ref{gasconstr}, but, for convenience, recall the definitions of the sets $Y_{\ad}$, that is, pointwise in time gradient constraints:
\[Y_{\ad} \coloneqq \left\{ y \in \W^{1,s}(I,H^{-\zeta,p}_D, \Dkal_\zeta)\colon\lVert \nabla y \rVert_{L^q}^q \leq g_{\avg} \text{ a.e.\ on } I \right\}\tag{\ref{Yadav}}\]
and fully pointwise ones:
\[ Y_{\ad} \coloneqq \left\{ y \in \W^{1,s}(I,H^{-\zeta,p}_D, \Dkal_\zeta)\colon \lvert \nabla y \rvert_2^2 \leq g \text{ a.e.\ on } Q \right\}.\tag{\ref{Yadptw}}\]

Let us now define the optimal control problem~\eqref{ocp} formally in reduced form, including the implicit constraint defined by $U_g$. By slight abuse of notation, we still refer to this problem as~\eqref{ocp}, taking for granted that the following is the---for the present purpose---correct and in particular mathematically well defined interpretation. From Proposition~\ref{prop::openset}, recall that the control-to-state operator $S$ is well defined and continuously differentiable on an open subset $U_g$ of $L^s(\Lambda)$. This allows to introduce the reduced cost functional
\[
j\colon U_g \rightarrow \R, \qquad u \mapsto J(S(u),u),
\]
that is also continuously differentiable due to the chain rule,
and the reduced control problem formulation
\[
\min_{u} j(u)\qquad \text{s.t.} \quad u\in U_{\ad} \cap U_g,\quad S(u)\in Y_{\ad}.\tag{\ref{ocp}}
\]
Again, note that $U_g$ is an \emph{open} subset of $L^s(\Lambda)$, that is, it is a continuum and sufficiently rich to do optimization on---provided that it is nonempty, and compatible with $U_{\ad}$ and $Y_{\ad}$. Regarding this  issue, we further pose the following usual assumption that, roughly speaking, ensures compatibility of the various constraints, enforcing a nonempty feasible set for~\eqref{ocp}.
\begin{assumption}\label{gasfeasibility}
Assume that there exists a feasible control-state-pair for~\eqref{ocp}, that is, there exists $\hat u \in U_{\ad} \cap U_g$ such that $S(\hat u) \in Y_{\ad}$.
\end{assumption}
In general, such an assumption cannot be avoided in PDE-constrained optimization with additional constraints on the state---see for example~\cite[Theorem~3]{CasasFernandez1993}---and in fact we will require even a bit more for necessary optimality conditions in Section~\ref{sec::fon} below, relying on Slater type conditions. For some exceptional special cases, in which a feasible control-state-pair and/or a Slater point can be constructed explicitly, we refer the reader to, e.g.,~\cite[Theorem~2.1]{Casas2014} for the case of a linear elliptic PDE\@. Due to the rather involved structure of the state equation in this work we only give
the following simple example for the model nonlinearity~\eqref{eq:examplesF1}:

\begin{example}
Let $u_{\low} \leq 0 \leq u_{\up}$, $\Fkal(y) = \lvert \nabla y \rvert^2$, and $y_0 \equiv 0$. Then, $\hat y \equiv 0$ is a global-in-time solution associated to the control $\hat u \equiv 0$, and it follows that $(\hat y,\hat u) \in Y_{\ad} \times U_{\ad}$. Consequently, Assumption~\ref{gasfeasibility} is fulfilled.
\end{example}

\subsection{Existence of globally optimal solutions}

Next, we collect the essential closedness property with respect to the weak topology on $L^s(\Lambda)$ of the feasible set in~\eqref{ocp} that allows to prove existence of optimal controls.

\begin{lemma} Either set $Y_{\ad}$ defined in~\eqref{Yadav} or~\eqref{Yadptw} in Assumption~\ref{gasconstr}, respectively, is weakly closed in $\W^{1,s}(I,H^{-\zeta,p}_D,\Dkal_\zeta)$.\label{lem:yad-weakly-closed}\end{lemma}

\begin{proof}
Due to compactness of the embedding~\eqref{eq:maxreg-into-cw1q}, it is sufficient to establish that either $Y_{\ad}$ is closed with respect to the $C^\alpha(\overline I,W^{1,q}_D)$-topology. For~\eqref{Yadav}, this is clear, whereas for~\eqref{Yadptw},
we note that if a sequence $(y_k)$ is convergent in $C^\alpha(\overline I,W^{1,q}_D)$ with limit $\bar y$, then the sequence of squared gradients $(\abs{\nabla y_k}^2)$ is convergent in $L^{q/2}(Q)$ with the limit $\abs{\nabla \bar y}^2$. It follows that there is a subsequence of $(\abs{\nabla y_k}^2)$ that converges pointwise almost everywhere on $Q$ to $\abs{\nabla \bar y}^2$. In particular, if $g$ was a uniform upper bound for each $\abs{\nabla y_k}^2$, then so is it for $\abs{\nabla \bar y}^2$, so $Y_{\ad}$ is indeed closed with respect to the $C^\alpha(\overline I,W^{1,q}_D)$-topology.
\end{proof}

\begin{lemma}\label{lem::closednessadm}
The feasible set of~\eqref{ocp}, that is, the set of all $u \in U_{\ad} \cap U_g$ such that $S(u) \in Y_{\ad}$, is weakly closed in $L^s(\Lambda)$.
\end{lemma}

\begin{proof}
 Let $(u_k) \subseteq U_{\ad} \cap U_g$ be a weakly convergent sequence in $L^s(\Lambda)$ such that $(y_k) \subseteq Y_{\ad}$, with $y_k \coloneqq S(u_k)$, and denote its limit by $\bar u$. Clearly, $U_{\ad}$ is convex and closed in $L^s(\Lambda)$ and thus weakly closed, so already $\bar u \in U_{\ad}$. In order to show that $\bar u \in U_g$ and $\bar y = S(\bar u) \in Y_{\ad}$, we rely on the interplay of $(y_k) \subseteq Y_{\ad}$ and Assumption~\ref{gassemilin}~\ref{gassemilin:boundedness} via Proposition~\ref{prop::closedness}. Indeed, Assumption~\ref{gassemilin}~\ref{gassemilin:boundedness} and the uniform gradient bounds on $(y_k)$ implied by $(y_k) \subseteq Y_{\ad}$ yield that the sequence $(\Fkal(y_k))$ is bounded in $L^s(I,H^{-\zeta,p}_D)$. But then Proposition~\ref{prop::closedness} strikes to show that $\bar u \in U_g$ and $y_k \rightharpoonup \bar y = S(\bar u)$ in $\W^{1,s}(I,H^{-\zeta,p}_D,\Dkal_\zeta)$, and Lemma~\ref{lem:yad-weakly-closed} finishes the proof.
 \end{proof}

Our first main result
follows now easily.

\begin{theorem}\label{thm::existence}
Let the feasibility Assumption~\ref{gasfeasibility} hold true. Then, there exists at least one solution $\bar u \in U_{\ad}\cap U_g$ to~\eqref{ocp} with associated state $\bar y \in Y_{\ad}$.
\end{theorem}

\begin{proof}
  With Assumption~\ref{gasfeasibility} and Lemma~\ref{lem::closednessadm} the result follows by standard arguments.
\end{proof}

To conclude this section we note that Theorem~\ref{thm::existence} stays true if $U_{\ad}$ is unbounded in $L^s(\Lambda)$, but in return the objective function is $L^s(\Lambda)$-coercive, for example by using a suitable regularization in $J$; in the context of quasilinear parabolic state equations, see e.g.~\cite[Proposition~6.4]{Meinlschmidt2016} or~\cite[Section~6]{HoppeNeitzel2020} for this approach.

\section{First-order optimality conditions}\label{sec::fon}

We now turn to the analysis of first-order necessary optimality conditions (FONs) characterizing locally
optimal solutions of~\eqref{ocp}. To this end, we recall the following standard notion, tacitly using the
notion $\B_\varepsilon(\bar u)$ for the open ball of radius $\varepsilon$ around $\bar u$ in
$L^s(\Lambda)$, here and in all that follows.

\begin{definition}We call $\bar u$ a locally optimal solution of~\eqref{ocp} if it is feasible, that is, $\bar u \in U_g \cap U_{\ad}$ and $S(\bar u) \in Y_{\ad}$, and there is some $\varepsilon > 0$ such that $j(u) \geq j(\bar u)$ holds for all feasible $u \in \B_\varepsilon(\bar u)$.\end{definition}
We still rely on the rather general regularity assumptions from Section~\ref{sec:notass} for the problem data and the choice of $p$ and $q$, as well as possible associated restrictions on $\zeta$ from the beginning of Section~\ref{sec::existence}, respectively, Assumption~\ref{gassemilin}~\ref{gassemilin:boundedness}. Again, these assumptions facilitate that
\begin{equation*}
   \W^{1,s}(I,H^{-\zeta,p}_D,\Dkal_\zeta) \hookrightarrow_c C^\alpha\bigl(\overline I,W^{1,q}_D\bigr) \tag{\ref{eq:maxreg-into-cw1q}}
\end{equation*}
precisely for the integrability index $q$ set forth by Assumption~\ref{gassemilin}~\ref{gassemilin:boundedness} which is thus also present in the gradient constraints as in Assumption~\ref{gasconstr}. Let us also recall that in the irregular case \irreg{} of Assumption~\ref{assu:regularity-settings}, we require $q=p$. Further, as one may expect in the context of deriving first order optimality conditions, we will now additionally rely on Assumption~\ref{gassemilin}~\ref{gassemilin:diffbar} that ensures differentiability of $\Fkal$. That is, we suppose that \emph{all} assumptions in Section~\ref{sec:notass} hold true for this section.

In the following, we discuss FONs for both possible choices to realize gradient constraints as provided by $Y_{\ad}$ in Assumption~\ref{gasconstr}, that is, pointwise in space and time,~\eqref{Yadptw}, as well as integral in space and pointwise in time,~\eqref{Yadav}. A particular aspect to consider is the nonstandard constraint set $U_g$ used to make~\eqref{ocp} well defined. However, from Proposition~\ref{prop::openset}, we do already know that under Assumption~\ref{gassemilin}~\ref{gassemilin:diffbar} $U_g$ is \emph{open}, so it is \emph{locally never active} in a local optimal solution to~\eqref{ocp} and will thus be irrelevant for optimality conditions. Mathematically, we can get rid of it utilizing a localization technique which was employed similarly before in~\cite{Meinlschmidt2017_2}; see Remark~\ref{rem:slater-for-global-control} below.

We give FONs separately for both choices of $Y_{\ad}$ in Assumption~\ref{gasconstr} as they yield structurally different optimality conditions. To expand on this further, it is well known that for optimal control problems involving pointwise constraints on the state, in order to employ a linearized Slater type constraint qualification, one needs, roughly speaking, the state to be continuous in any component that is constrained in a pointwise manner, and the associated Lagrangian multiplier will be a measure in that component. In the present case, in the case of integral in space but pointwise in time gradient constraints as in~\eqref{Yadav}, the gradient of the state would have to be continuous in the time component with values in $L^q$---see~\eqref{eq:maxreg-into-cw1q}---and the Lagrange multiplier would be a measure in time, whereas in case of full pointwise constraints on $Q$ on the gradient of the state as in~\eqref{Yadptw}, we would require continuity of the gradient of the state on $Q$ and obtain a Lagrange multiplier which is represented by a measure on $Q$. There seems to be no obvious way of stating a general flexible result which works for both cases, so it makes sense to branch here. In both cases, we will eventually also slightly sharpen the assumptions on $g_{\avg}$, $g$, respectively, assuming that they are continuous and strictly positive functions; cf.\ Definition~\ref{def:linearized-slater} and Remark~\ref{rem:linearized-slater-inactive-zero}. Note that this does not affect the previous results on existence of solutions if Assumption~\ref{gasfeasibility} holds accordingly.

We continue with the linearized Slater constraint qualifications. It will be useful to abbreviate the following (nonlinear) operator, the $q$-Laplacian with factor $q$,
\begin{equation}
  -\Delta_q \colon W^{1,q}_D \to W^{-1,q'}_D, \qquad \langle -\Delta_q \varphi, \psi\rangle \coloneqq q\int_\Omega \lvert \nabla \varphi \rvert^{q-2}\nabla \varphi \cdot \nabla \psi \d x,\label{eq:q-laplace}
\end{equation}
which arises in the present context as the derivative of $W^{1,q}_D \ni z \mapsto \lVert \nabla z \rVert_{L^q}^q$ defined in~\eqref{eq:grad-lq-norm}.

\begin{definition}[Linearized Slater condition]\label{def:linearized-slater} Let $\bar u$ be a feasible point for~\eqref{ocp} with associated state $\bar y = S(\bar u)$. Then we say that the \emph{linearized Slater condition} is satisfied at $\bar u$ if there exists $\hat u \in U_{\ad}$ and $\tau > 0$ such that:
\begin{enumerate}[(i)] \item\label{def:linearized-slater-1}
If $Y_{\ad}$ is given by~\eqref{Yadav} with $g_{\avg} \in C(\overline I)$,
\begin{equation*}
\lVert \nabla \bar y (t) \rVert_{L^q}^q +
\bigl\langle-\Delta _q \bar y(t), z(t)\bigr\rangle
\leq g_{\avg}(t) - \tau \qquad \text{for all}~t \in \overline I,
\end{equation*}
\item\label{def:linearized-slater-2} if $Y_{\ad}$ is given by~\eqref{Yadptw} with $g \in C(\overline Q)$,
\begin{equation*}
\lvert \nabla \bar y \rvert^2 +
2 \nabla \bar y \cdot \nabla z
\leq g - \tau \qquad \text{on}~\overline Q,
\end{equation*}
\end{enumerate}
where in either case $z = S'(\bar u)(\hat u-\bar u) \in \W^{1,s}(I,H^{-\zeta,p}_D,\Dkal_\zeta)$ is the unique solution to~\eqref{lse} with $v = \hat u-\bar u$, that is, to
\begin{align*}
\begin{aligned}
\partial_t z + \Akal(\bar y) z + \Alin(\bar y)z -\Fkal'(\bar y) z &= B (\hat u - \bar u) && \text{in } H^{-\zeta,p}_D \text{ a.e.\ on } I,\\
z(0) &= 0 && \text{in } (H^{-\zeta,p}_D,\Dkal_\zeta)_{1-\frac1s,s}.
\end{aligned}
\end{align*}
\end{definition}

\begin{remark}\label{rem:linearized-slater-inactive-zero}
  The linearized Slater condition contains several implicit conditions:
  \begin{enumerate}[(i)]
  \item Clearly, in order for the condition to hold, we must necessarily have $g_{\avg} > 0$ uniformly on $\overline I$ and $g > 0$ uniformly on $\overline Q$, respectively. Since we have included uniform continuity of either function in the linearized Slater condition, this means they will have a strictly positive lower bound.
   \item Since $z(0) = 0$, the linearized Slater condition requires that the respective constraint is \emph{inactive} for $t = 0$. Since $\bar y(0) = y_0$ with $\bar y = S(\bar u)$, that means that necessarily $\norm{\nabla y_0}_{L^q}^q < g_{\avg}(0)$ when $Y_{\ad}$ is given by~\eqref{Yadav} and $\abs{\nabla y_0}^2 < g(0)$ uniformly on $\overline\Omega$ when $Y_{\ad}$ is given by~\eqref{Yadptw}.
  \end{enumerate}
\end{remark}

\begin{remark}\label{rem:slater-for-global-control}
We remark that $U_g$ does \emph{not} occur in Definition~\ref{def:linearized-slater} at all. In fact, it is easily seen that if the linearized Slater condition as in Definition~\ref{def:linearized-slater} is satisfied in $\bar u$ which is feasible for~\eqref{ocp}, then, for any $\delta > 0$, we can without loss of generality choose an associated linearized Slater point which is in $U_{\ad} \cap \mathbb{B}_\delta(\bar u)$. Since $\bar u \in U_g$ and the latter is open, this implies in particular that we can pick a linearized Slater point in $U_{\ad} \cap U_g$. To confirm this, suppose that $\hat u \in U_{\ad}$ is the linearized Slater point as in Definition~\ref{def:linearized-slater} with associated $z = S'(\bar u)(\hat u - \bar u)$ and set $\hat u_\alpha \coloneqq (1-\alpha) \bar u + \alpha \hat u$ for $\alpha \in (0,1)$. Observe that $\hat u_\alpha \in U_{\ad}$ since $U_{\ad}$ was convex, and that $\hat u_\alpha - \bar u = \alpha(\hat u-\bar u)$, so that $z_\alpha \coloneqq S'(\bar u)(\hat u_\alpha - \bar u) = \alpha z$. In particular, for a fixed $\delta > 0$ we find $\alpha = \alpha(\delta) > 0$ sufficiently small such that $\hat u_\alpha \in U_{\ad} \cap \mathbb{B}_\delta(\bar u)$ and further
\begin{equation*}
\lVert \nabla y (t) \rVert_{L^q}^q +
\bigl\langle-\Delta _q y(t), z_\alpha(t)\bigr\rangle
\leq g_{\avg}(t) - \alpha\tau \qquad \text{for all}~t \in \overline I
\end{equation*}
in the case~\eqref{Yadav} for $Y_{\ad}$ or, respectively, in the case~\eqref{Yadptw} for $Y_{\ad}$,
\begin{equation*}
\lvert \nabla \bar y \rvert^2 +
2 \nabla \bar y \cdot \nabla z_\alpha
\leq g - \alpha\tau \qquad \text{on}~\overline Q.
\end{equation*}
Hence, $\hat u_\alpha \in U_{\ad} \cap \mathbb{B}_\delta(\bar u)$ also provides the linearized Slater condition for the feasible set of~\eqref{ocp} as in Definition~\ref{def:linearized-slater}. Since $\bar u$ is feasible for~\eqref{ocp} and $U_g$ is open by Proposition~\ref{prop::openset}, for $\delta > 0$ sufficiently small, we have $\mathbb{B}_\delta(\bar u) \subseteq U_g$ and thus also $\hat u_\alpha \in U_{\ad} \cap U_g$.
\end{remark}

Before we continue with our second main result, we briefly clarify the notion of very weak solutions of (adjoint) equations.
For that, recall the solution operator for the linearized state equation \begin{equation*}
  \altS(y) \coloneqq \bigl(\partial_t + \Akal(y) + \Alin(y) - \Fkal'(y)\bigr)^{-1} \colon L^s(I,H^{-\zeta,p}_D)
  \to \W^{1,s}_0(I,H^{-\zeta,p}_D, \Dkal_\zeta).
\end{equation*} as in Remark~\ref{rem:general-linearized-solution-operator} for $y \in \W^{1,s}(I,H^{-\zeta,p}_D,\Dkal_\zeta)$.

\begin{definition}[Very weak solution, adjoint equation]\label{def:very-weak-adjoint}
Suppose that $X$ is a Banach space such that $ \W^{1,s}_0(I,H^{-\zeta,p}_D,\Dkal_\zeta) \hookrightarrow C(\overline I,X)$ and denote that embedding by $\Ekal$.
Let $y \in \W^{1,s}(I,H^{-\zeta,p}_D,\Dkal_\zeta)$ and $\eta \in \Mkal(\overline I,X^*)$ be given, the latter in the form $\eta = \eta_I + \delta_T \otimes \eta_T$ with $\eta_T \in X^*$ and $\eta_I(\{T\}) = 0$. Then we say that $w \in L^{s'}(I,H^{\zeta,p'}_D)$ is the \emph{very weak} solution to the \emph{adjoint equation}
\begin{equation}\left.\begin{aligned}
  -\partial_t w + \Akal(y)^*w + \Alin(y)^* w - \Fkal'(y)^* w & = \eta_I, \\ w(T) & = \eta_T
\end{aligned}\quad \right\}\label{eq:general-adjoint}
  \end{equation}
  if $w = \altS(y)^*\Ekal^*\eta$, i.e., if, for all $v \in \W^{1,s}(I,H^{-\zeta,p}_D,\Dkal_\zeta)$,
   \begin{equation*}
\int_0^T \bigr\langle (\partial_t + \Akal(y) + \Alin(y) - \Fkal'(y))v,w\bigr\rangle_{H^{-\zeta,p}_D,H^{\zeta,p'}_D} \\ = \int_0^T v \d \eta_I + \langle \eta_T,v(T)\rangle_{X^*,X}.
  \end{equation*}
\end{definition}

\begin{remark}
Clearly, some comments on Definition~\ref{def:very-weak-adjoint} are appropriate.
\begin{enumerate}[(i)]
  \item First, $w$ is well defined. Indeed, we have the topological isomorphism
  \begin{equation*}\altS(y)^* \colon
  \W^{1,s}_0(I,H^{-\zeta,p}_D,\Dkal_\zeta)^* \to L^s(I,H^{-\zeta,p}_D)^* = L^{s'}(I,H^{\zeta,p'}_D)
  \end{equation*} and the adjoint of the embedding
  \[\Ekal^* \colon \Mkal(\overline I,X^*) = C(\overline I,X)^* \hookrightarrow \W^{1,s}_0(I,H^{-\zeta,p}_D,\Dkal_\zeta)^*,\]
  so that $\altS(y)^*\Ekal^*\eta \in L^{s'}(I,H^{\zeta,p'}_D)$. Since $\altS(y)^*$ is invertible, the very weak solution is also in fact unique.
  \item We hope that the reader agrees with our assessment of a \emph{very weak} sense of solution. It is well known that one cannot expect a weak time derivative for $w$ in the present context with measure data. Now, \emph{if} $\eta_I \in L^s(I,\Dkal_\zeta^*)$ and $(\Akal(y) + \Alin(y) - \Fkal'(y))^*$ satisfies nonautonomous maximal parabolic regularity on $\Dkal_\zeta^*$ with constant domain $H^{\zeta,p'}_D$, \emph{then} we in fact obtain that $w \in \W^{1,s'}(I,\Dkal_\zeta^*,H^{\zeta,p'}_D)$ and~\eqref{eq:general-adjoint} in the strong sense in $\Dkal_\zeta^*$; see~\cite[Section~6]{Amann2005}. We do not go into the details of the adjoint equation analysis here since in the first order necessary conditions below, we will be limited by the regularity of $\eta$, which will indeed be a measure only, anyhow.
\end{enumerate}
\end{remark}

Now let us finally turn to the actual optimality conditions for each case of $Y_{\ad}$.

\subsection{Integral in space and pointwise in time gradient-constraints}\label{sec:fon-integr}

We first consider the case of integral in space and pointwise in time gradient-constraints as in~\eqref{Yadav} in Assumption~\ref{gasconstr}, that is,
\[Y_{\ad} \coloneqq \left\{ y \in \W^{1,s}(I,H^{-\zeta,p}_D, \Dkal_\zeta)\colon\lVert \nabla y \rVert_{L^q}^q \leq g_{\avg} \text{ a.e.\ on } I \right\}.\]
In accordance with the convention at the end of Assumption~\ref{gasconstr}, we tacitly suppose $q<\infty$ for this section; the case $q=\infty$ which effectively results in pointwise in time and space constraints as in~\eqref{Yadptw} will be covered in Section~\ref{sec::higher} below.

The following first-order necessary optimality conditions are our second main result.
Let us emphasize that we can afford the irregular setting \irreg{} of Assumption~\ref{assu:regularity-settings} with $\frac{d}p < \zeta < 1-\frac2s$ and $q=p$ for this type of constraint.

A similar result has been obtained in~\cite[Theorem~2.5]{Ludovici2015} for linear parabolic problems and $p=2$ with pointwise-in-time but integrated-in-space gradient constraints; see also the extension to semilinear problems in~\cite{LudoviciNeitzelWollner}.

\begin{theorem}[First order optimality conditions/KKT conditions]\label{thm::fonintegr}
Consider the case of Assumption~\ref{gasconstr}~\ref{gasconstr-avg} on $Y_{\ad}$ with $0 < g_{\avg} \in C(\overline I)$. Let $\bar u \in U_{\ad} \cap U_g$ be a local solution to~\eqref{ocp} with associated state $\bar y = S(\bar u)$ for which the linearized Slater condition as in Definition~\ref{def:linearized-slater}~\ref{def:linearized-slater-2} is satisfied.

Then, there exists a regular Borel measure $\bar \lambda \in \Mkal(\overline I) = C(\overline I)^*$, the Lagrangian multiplier associated to the gradient constraint, such that the \emph{complementarity condition}
\begin{equation}\label{fonintegr_compl}
\int_I \bigl( \varphi - \lVert \nabla \bar y \rVert_{L^q}^q \bigr) \d\bar \lambda \leq 0 \qquad \text{for all } \varphi \in C(\overline I) \text{ such that } \varphi \leq g_{\avg}\end{equation}
and the \emph{variational inequality} are satisfied:
\begin{equation}\label{fonintegr_vi}\int_\Lambda (\gamma \bar u + B^* \bar w)(u-\bar u) \d \rho
\geq 0 \qquad \text{for all } u \in U_{\ad}.
\end{equation}
Here, $\bar w \in L^{s'}(I,H^{\zeta,p'}_D)$ is the \emph{adjoint state} defined as the very weak solution of  the abstract \emph{adjoint equation} in the sense of Definition~\ref{def:very-weak-adjoint},
\begin{equation}\label{fonintegr_ae} \left. \begin{aligned}
-\partial_t \bar w + \Akal(\bar y)^* \bar w + \Alin(\bar y)^* \bar w - \Fkal'(\bar y)^* \bar w &= \bar y - y_d
  + \d\bar\lambda_I \otimes (-\Delta_q \bar y), \\
w(T) &= -\bar\alpha_T\Delta_q \bar y(T),
\end{aligned}\quad \right\}
\end{equation}
where we have decomposed $\bar \lambda = \bar\lambda_I + \bar\alpha_T \delta_T$ with $\bar \lambda_I(\{T\}) = 0$ and $\bar \alpha_T \in \R$.
\end{theorem}

\begin{proof}
  The proof of Theorem~\ref{thm::fonintegr} is accomplished by an application of a well known abstract result for problems with state constraints due to Casas~\cite[Theorem~5.2]{Casas1993}. However, the constraint set for controls $U_{\ad} \cap U_g$ of~\eqref{ocp} is in general nonconvex, which does not fit with the assumptions of the mentioned theorem and general KKT theory. To rectify this issue, we follow the ideas of~\cite[Theorem 4.8]{Meinlschmidt2017_2} and localize the problem: Since $U_g$ is open in $L^s(\Lambda)$ according to Proposition~\ref{prop::openset} and $\bar u$ is a local solution of~\eqref{ocp}, we find  $\varepsilon > 0$ sufficiently small such that $\mathbb{B}_\varepsilon(\bar u) \subset U_g$ and such that $\bar u$ is a solution of~\eqref{ocp} with $U_{\ad}$ replaced by $U_{\ad} \cap \mathbb{B}_\varepsilon(\bar u)$. The latter is a \emph{convex} set, and via Remark~\ref{rem:slater-for-global-control} we see that the linearized Slater constraint qualification at $\bar u$ is still satisfied for the localized problem. Thus, we can apply~\cite[Theorem~5.2]{Casas1993} to the latter by setting $U = L^s(\Lambda)$, $Z = C(\overline I)$, $K = U_{\ad} \cap \mathbb{B}_\varepsilon(\bar u)$, $C \coloneqq \{ \varphi \in Z\colon \varphi \leq g_{\avg} \text{ on } \overline I \}$, $J = j$, and $G = \mathcal{I} \circ \Ekal \circ S$. Hereby, $j$ denotes the reduced functional of~\eqref{ocp}, which is continuously differentiable on $U = L^s(\Lambda)$, $\Ekal$ is the embedding $\W^{1,s}(I,H^{-\zeta,p}_D,\Dkal_\zeta) \hookrightarrow C(\overline I,W^{1,q}_D)$ as in~\eqref{eq:maxreg-into-cw1q}, and
\[
\mathcal{I}\colon C(\overline I,W^{1,q}_D) \rightarrow C(\overline I), \qquad y \mapsto \lVert \nabla y(\cdot) \rVert_{L^q}^q \coloneqq \int_{\Omega} \lvert \nabla y(\cdot) \rvert^q \d x,
\]
is taking spatial $L^q$-integrals of the gradients. We note that $\mathcal{I}$ is continuously Fréchet differentiable with the derivative in direction $z$ as an element of $C(\bar I)$ given by
$\mathcal{I}'(\bar y)z = \langle -\Delta_q y(\cdot),z(\cdot)\rangle$, with $-\Delta_q$ as in~\eqref{eq:q-laplace}.
It follows that $G$ is continuously differentiable with $G'(\bar u) = \mathcal{I}'(\bar y) \Ekal S'(\bar u)$. Note that $S'(\bar u) = \altS(\bar y) B$, see Remark~\ref{rem:general-linearized-solution-operator}, and thus $G'(\bar u) = \mathcal{I}'(\bar y) \Ekal\altS B $. 
From~\cite[Theorem 5.2]{Casas1993} we obtain a Lagrange multiplier $\bar \lambda \in \Mkal(\overline I) = C(\overline I)^*$ such that the complementarity condition~\eqref{fonintegr_compl} is satisfied, as well as the variational inequality
\begin{equation*}
\int_\Lambda (\gamma \bar u + B^* \bar w)(u-\bar u) \d \rho
\geq 0 \qquad \text{for all } u \in U_{\ad} \cap \mathbb{B}_\varepsilon(\bar u),
\end{equation*}where $\bar w$ is the adjoint state
\begin{align}\label{defw}
\bar w \coloneqq \altS(\bar y)^*\Ekal^* \Bigl[\bar y-y_d + \mathcal{I}'(\bar y)^* \bar \lambda\Bigr] \in L^{s'}(I,H^{\zeta,p'}_D).
\end{align}
From the foregoing localized variational inequality, we recover the global one in~\eqref{fonintegr_vi} which is valid for all $u \in U_{\ad}$ by a simple rescaling argument using convexity of $U_{\ad}$. Note that adjoint state in~\eqref{defw} is well defined, but we have tacitly identified $\bar y - y_d \in L^\infty(I,L^2)$ as an element of $C(\overline I,W^{1,q}_D)^*$.

It remains to dissect $\Ikal'(\bar y)^*\bar\lambda \in C(\overline I,W^{1,q}_D)^* = \Mkal(\overline I,W^{-1,q'}_D)$. From the linearized Slater condition for $\bar u$ via Remark~\ref{rem:linearized-slater-inactive-zero}, and continuity of the functions involved, it follows that there is $\eps > 0$ such that the gradient constraint is inactive on $[0,2\eps)$. Using the complementarity condition~\eqref{fonintegr_compl}, we easily infer that $\supp \bar\lambda \cap [0,\eps] = \emptyset$, in particular, $0 \notin \supp \bar\lambda$. We obtain the decomposition $\bar \lambda = \bar\lambda_I + \bar\alpha_T \delta_T$ by restriction of $\bar\lambda$ to $I$ and $\{T\}$, respectively, with $\supp \bar\lambda_I \subset I$ and $\bar\alpha_T \in \R$. (In fact, $\bar\lambda$ is a nonnegative measure and $\bar\alpha_T \geq 0$; see Remark~\ref{rem:fonintegr} below.) Now, turning to $\Ikal'(\bar y)^*\bar\lambda$, for any $z \in C(\overline I, W^{1,q}_D)$,
\begin{align*}
\bigl\langle \mathcal{I}'(\bar y)^* \bar \lambda, z \bigr\rangle = \bigl\langle \bar\lambda, \Ikal(\bar y)'z\bigr\rangle
= \int_{I} \langle -\Delta_q \bar y,z\rangle \d\bar\lambda_I - \bar\alpha_T\bigl\langle\Delta_q \bar y(T),z(T)\bigr\rangle,
\end{align*}
which we denote by $\mathcal{I}'(\bar y)^* \bar\lambda =
\d\bar\lambda_I \otimes (-\Delta_q \bar y) + \bar\alpha_T \delta_T \otimes (-\Delta_q\bar y(T))$. Doing so, it follows that $\bar w$ is a very weak solution to the adjoint equation~\eqref{fonintegr_ae} as in Definition~\ref{def:very-weak-adjoint} with $X = W^{1,q}_D$.
\end{proof}

\begin{remark}\label{rem:fonintegr} Let us state some rather typical conclusions from Theorem~\ref{thm::fonintegr}.
\begin{enumerate}[(i)]
\item As usual, the variational inequality~\eqref{fonintegr_vi} can be expressed equivalently as \[\bar u(z) = \operatorname{proj}_{[u_{\low}(z),u_{\up}(z)]}\bigl(-\gamma^{-1}(B^*\bar w)(z)\bigr) \qquad \text{for $\rho$-almost all $z \in \Lambda$}.\]
\item The complementarity condition~\eqref{fonintegr_compl} implies that $\bar\lambda \geq 0$ in the sense of measures in $\Mkal(\overline I)$, that is, $\langle \bar\lambda,\psi\rangle \geq 0$ for all $\psi \in C(\overline I)$ with $\psi \geq 0$. Choose $\varphi \coloneqq \norm{\nabla \bar y}_{L^q}^q - \psi$ in~\eqref{fonintegr_compl} to see this. This also means that both $\bar\lambda_I$ and $\bar\alpha_T$ are also nonnegative. In fact, it is easily seen that the complementarity condition~\eqref{fonintegr_compl} is equivalent to the classical complementarity condition
\begin{equation*}
  \bar\lambda \geq 0 \quad \text{and} \quad \int_{\bar I} \bigl(g_{\avg} - \norm{\nabla \bar y}_{L^q}^q\bigr) \d\bar\lambda = 0.
\end{equation*}
In particular, it follows that \[\supp \bar \lambda \subseteq \Bigl\{t \in \overline I \colon \norm{\nabla \bar y(t)}_{L^q}^q = g_{\avg}(t)\Bigr\}.\]
\end{enumerate}
\end{remark}

\subsection{Pointwise in space and time constraints on the gradient}\label{sec::higher}
In this section we derive first-order optimality conditions for~\eqref{ocp} with pointwise in space and time gradient-constraints as in~\eqref{Yadptw} in Assumption~\ref{gasconstr}, that is,
\[ Y_{\ad} \coloneqq \left\{ y \in \W^{1,s}(I,H^{-\zeta,p}_D, \Dkal_\zeta)\colon \lvert \nabla y \rvert_2^2 \leq g \text{ a.e.\ on } Q \right\}.\]
As already explained at the beginning of Section~\ref{sec::fon}, the linearized Slater condition as in Definition~\ref{def:linearized-slater} requires to consider the states $y = S(u)$ in a function space whose elements have gradients which are continuous on $\overline Q$. This will in general not be achievable within the irregular framework \irreg{} in Assumption~\ref{assu:regularity-settings}. Consequently, in this section, in contrast to Section~\ref{sec:fon-integr} and irrespective of the value of $q$ in Assumption~\ref{gassemilin}~\ref{gassemilin:boundedness}, we restrict ourselves to the regular case \reg{} in Assumption~\ref{assu:regularity-settings} and $\zeta < 1-\frac{d}p-\frac{2}s$, from which we can leverage the required regularity
\begin{equation*}
   \W^{1,s}(I,H^{-\zeta,p}_D,\Dkal_\zeta) \hookrightarrow_c C^\alpha\bigl(\overline I,C^1(\overline\Omega)\bigr).
\end{equation*}
by means of Corollary~\ref{cor:MPR-cont-W1q-embed}. We point out that this requirement fits precisely with the interpretation of the gradient constraint sets $Y_{\ad}$ in Assumption~\ref{gasconstr} for $q=\infty$; in particular, it works for every $q$ given by Assumption~\ref{gassemilin}~\ref{gassemilin:boundedness}, including $q=\infty$. In summary:
\begin{center}
  \vspace{0.75em}
  \fbox{\parbox{0.95\linewidth}{
  For this section, irrespective of the value of $q$ in Assumption~\ref{gassemilin}~\ref{gassemilin:boundedness}, we require the regular setting \reg{} with \[\domLp = W^{2,p} \cap W^{1,p}_D \qquad \text{and} \qquad \zeta < 1-\frac2s - \frac{d}p.\]
  }}
  \vspace{0.75em}
\end{center}

\begin{remark}\label{rem-comments-on-irreg}
Let us mention again that, as outlined in Example~\ref{ex:data-regularity-strong}, the implicit implications of the regular setting \reg{} are not insignificant as they, in general, exclude mixed boundary conditions that meet, nonsmooth domains, and non-Lipschitz coefficients. From the examples for $B$ summarized below Assumption~\ref{gaseq}, cf.\ also in~\cite{Bonifacius2018,HoppeNeitzel2020}, the following analysis is effectively restricted to distributed control on (subsets of) the space-time cylinder $Q$, or controls acting in time only with fixed actuators from, say, $H^{-\zeta,p}_D$. However, when discussing pointwise constraints on the gradient of the state, such a restriction to a rather regular setting cannot be avoided.
\end{remark}

With the regularity of the states provided by the regularity setting and the choice of $\zeta$, we can again apply~\cite[Theorem~5.2]{Casas1993} to~\eqref{ocp} with the present choice of pointwise constraints on the gradients of the states.
A result analogous to the following one, but for semilinear elliptic state equations, has been obtained in~\cite[Corollary~1]{CasasFernandez1993}.

\begin{theorem}\label{thm::fonptw}
  Consider the case of Assumption~\ref{gasconstr}~\ref{gasconstr-ptw} on $Y_{\ad}$ with $0 < g\in C(\overline Q)$. Let $\bar u \in U_{\ad} \cap U_g$ be a local solution to~\eqref{ocp} with associated state $\bar y = S(\bar u)$ for which the linearized Slater condition as in Definition~\ref{def:linearized-slater}~\ref{def:linearized-slater-1} is satisfied. Then, there exists a regular Borel measure $\bar \nu \in \Mkal(\overline Q) = C(\overline Q)^*$, the {Lagrangian multiplier} associated with the gradient constraint, such that the \emph{complementarity condition}
\begin{align}\label{fonptw_compl1}
\int_{\overline Q}\bigl(\varphi - \abs{\nabla y}_2^2\bigr) \d\bar\nu \leq 0, \qquad \text{for all } \varphi \in C(\overline Q)~\text{such that}~\varphi \leq g,
\end{align}
and the \emph{variational inequality} are satisfied:
\begin{align}\label{fonptw_vi}
\int_\Lambda (\gamma \bar u + B^* \bar w)(u-\bar u\bigr) \d \rho
\geq 0, \qquad \text{for all } u \in U_{\ad}.
\end{align}
Here, $\bar w \in L^{s'}(I,H^{\zeta,p'}_D)$ is the \emph{adjoint state} defined as the very weak solution of
the abstract \emph{adjoint equation} in the sense of Definition~\ref{def:very-weak-adjoint},
\begin{equation}\label{fonptw_ae}
\left. \begin{aligned}
-\partial_t \bar w + \Akal(\bar y)^* \bar w + \Alin(\bar y)^* \bar w - \Fkal'(\bar y)^* \bar w &= \bar y - y_d - 2 \nabla \cdot \bar \nu_Q \nabla \bar y, \\
\bar w(T) &= -2\nabla \cdot \bar\nu_T \nabla \bar y(T), \end{aligned} \quad \right\}
\end{equation}
where we have decomposed $\bar \nu = \bar \nu_Q + \delta_T \otimes \bar \nu_T$ with $\bar \nu_Q(\{T\}\times\overline\Omega) = 0$ and $\bar \nu_T \in \Mkal(\overline\Omega)$.
\end{theorem}

\begin{proof}
The proof is similiar to the proof of Theorem~\ref{thm::fonintegr}; see also~\cite[Theorem~4 and Corollary~1]{CasasFernandez1993}. We thus do not reiterate the details of the localization procedure---that is required for this case equally well---and restrict ourselves to the essential differences. This time, we apply~\cite[Theorem~5.2]{Casas1993} with the following choices, for $\eps>0$ sufficiently small: $Z = C(\overline Q)$, $U = L^s(\Lambda)$, $K=U_{\ad} \cap \B_\varepsilon(\bar u)$, $C= \{ \varphi \in C(\overline Q)\colon \varphi \leq g \}$, $J = j$, and $G = Q \circ \Ekal \circ S$ with $\Ekal$ being the embedding $\W^{1,s}(I,H^{-\zeta,p}_D,\Dkal_\zeta) \hookrightarrow C(\overline I,C^1(\overline\Omega))$ and $Q \colon z \mapsto \lvert \nabla z \rvert_2^2$ taking $C(\overline I,C^1(\overline\Omega) )$ into $C(\overline Q)$. Note that $\Ekal$ and $G$ are well defined due to Corollary~\ref{cor:MPR-cont-W1q-embed} and the choice of $\zeta$ for this section; moreover, $Q$ is continuously Fr{\'e}chet differentiable with the derivative $Q'(\bar y)z = 2 \nabla \bar y \cdot \nabla z$ in direction $z$. We obtain a Lagrangian multiplier $\bar \nu \in \Mkal(\overline Q) = C(\overline Q)^*$ for which---modulo the localization procedure---the complementarity condition~\eqref{fonptw_compl1} and the variational inequality~\eqref{fonptw_vi} are satisfied for the adjoint state
\begin{align*}
\bar w \coloneqq \altS(\bar y)^*\Ekal^*\Bigl[\bar y-y_d + Q'(\bar y)^* \bar \nu\Bigr] \in L^{s'}(I,H^{\zeta,p'}_D).
\end{align*}
Again, we consider $Q'(\bar y)^* \bar \nu$ in more detail. By restriction, we decompose the measure $\bar \nu = \bar \nu_Q + \delta_T \otimes \bar \nu_T$ with $\supp \bar \nu_Q \subseteq I \times \overline\Omega$ and $\bar \nu_T \in \Mkal(\overline\Omega)$. As before, the linearized Slater condition implies that $\bar \nu(\{0\} \times \overline\Omega) = 0$. Overall we find, for every $\psi \in C(\overline I,C^1(\overline\Omega))$,
\begin{equation*}
\bigl\langle Q'(\bar y)^* \bar \nu, \psi \bigr\rangle
= \bigl\langle \bar \nu, Q'(\bar y) \psi \bigr\rangle
= \int_{(0,T) \times \overline\Omega} 2 \nabla \bar y \cdot \nabla \psi \d\bar\nu_Q + \int_{\overline \Omega} 2 \nabla \bar y(T) \cdot \nabla \psi(T) \d \bar\nu_T.
\end{equation*}
We understand this in a weak sense $Q'(\bar y)^* \bar \nu = - 2\nabla \cdot \bar \nu_Q \nabla \bar y + 2\delta_T\otimes (-\nabla \cdot \bar\nu_T \nabla \bar y(T))$, and the adjoint equation~\eqref{fonptw_ae} follows as in Definition~\ref{def:very-weak-adjoint} with $X = C^1(\overline\Omega)$.
\end{proof}

\begin{remark}
The analogous comments as in Remark~\ref{rem:fonintegr} apply to Theorem~\ref{thm::fonptw} as well. In particular, there is also the pointwise representation for the optimal control $\bar u$ from the variational inequality~\eqref{fonptw_vi}, and the complementarity condition~\eqref{fonptw_compl1} is equivalently expressed by
\begin{equation*}
\bar \nu \geq 0 \quad \text{and} \quad \int_{\overline Q} \bigl(g - \abs{\nabla \bar y}_2^2\bigr) \d\bar\nu = 0
\end{equation*}
and we have
\begin{align*}
 \supp \bar \nu \subseteq \Bigl\{ (t,x) \in \overline Q\colon \lvert \nabla \bar y (t,x)\rvert_2^2 = g(t,x) \Bigr\}.
  \end{align*}
\end{remark}

\section{Extensions and variations of the results}\label{sec:generalize}
Recalling the precise constraint definitions in Assumption~\ref{gasconstr} with either plain pointwise in space and time gradients or pointwise in time but averaged in space gradient constraints combined with \emph{the condition $\Gamma_D\neq\emptyset$} and \emph{usage of Euclidean norms of the gradients}, let us now discuss some possible generalizations of our results.

\subsection{Pure Neumann boundary conditions in Assumption~\ref{gasconstr}}\label{sec:gamma}
So far, we have assumed $\Gamma_D\neq\emptyset$ in all definitions of $Y_{\ad}$, i.e.\ homogeneous Dirichlet boundary conditions were supposed to hold on a part $\Gamma_D$ of $\partial\Omega$ with nonzero surface measure. This was used in particular in the proof of existence of an optimal solution in Theorem~\ref{thm::existence} since the prerequisites of Propositions~\ref{prop::globalex} and~\ref{prop::closedness} guaranteeing global-in-time solutions of the state equation and appropriate convergence of infimizing sequences could be fulfilled.

In the case $\Gamma_D = \emptyset$, i.e.\ in the presence of pure Neumann boundary conditions, $L^q$-bounds on the gradient imposed via $Y_{\ad}$ do no more imply ``full'' $W^{1,q}$-bounds since the Poincar\'e-Friedrichs inequality mentioned in Remark~\ref{rem:semibounds} does not hold in that case. Thus, either $Y_{\ad}$ must include ``full'' $W^{1,q}$-bounds (which will be discussed in Subsection~\ref{sec:additional-constraints} below), or the assumptions on $\Fkal$ must be strengthened by modifying Assumption~\ref{gassemilin}~\ref{gassemilin:boundedness} as follows: for every
$C\geq 0$ there is $M \geq 0$ such that for every $y \in \W^{1,s}(I,H^{-\zeta,p}_D,\Dkal_\zeta)$:
\begin{equation*}
  \sup_{t\in \overline I} \lVert \nabla y(t) \rVert_{L^q} \leq C \quad \implies \quad \norm{\Fkal(y)}_{L^s(I,H^{-\zeta,p}_D)} \leq M.
\end{equation*}
In other words: If we assume boundedness of the $W^{1,q}$-seminorms of admissible $y$'s, uniform boundedness properties of $\Fkal$ have to be stated w.r.t.\ this $W^{1,q}$-seminorm accordingly. Clearly, this modified assumption would still be verified by the quadratic gradient nonlinearity~\eqref{eq:examplesF1}.

\subsection{Additional zero-order-state constraints}\label{sec:additional-constraints}

Now we will consider additional zero-order state constraints, i.e.\ the set $Y_{\ad}$ is defined by one of the following possible choices:
\begin{itemize}
\item Pointwise in space and time gradient-state-constraints:
\begin{align}\label{Yadptw0}Y_{\ad}& \coloneqq \left\{ y \in \W^{1,s}(I,H^{-\zeta,p}_D, \Dkal_\zeta)\colon \lvert \nabla y \rvert_2^2 \leq g, \quad y_{\low} \leq y \leq y_{\up} \; \text{a.e.\ on } Q \right\}
\end{align}
with $g \in L^\infty(I,L^{p/2})$, $y_{\low}, y_{\up} \in L^\infty(Q)$, and $g \geq 0$, $y_{\low} \leq y_{\up}$ a.e.\ on $Q$,
\item Integral in space and pointwise in time gradient-state-constraints:
\begin{multline}
Y_{\ad} \coloneqq \Bigl\{ y \in \W^{1,s}(I,H^{-\zeta,p}_D, \Dkal_\zeta)\colon \lVert \nabla y(t) \rVert_{L^q}^q \leq g_{\avg}(t),\\ \lVert y(t) \rVert_{L^q}^q \leq y_{\avg}(t) \quad \text{ a.e.\ on } I \Bigr\}\label{Yadav0}
\end{multline}
with $g_{\avg}, y_{\avg} \in L^\infty(I)$, $g_{\avg},y_{\avg} \geq 0$ a.e.\ on $I$.
\end{itemize}
This allows to consider both $\Gamma_D=\emptyset$ and $\Gamma_D\neq \emptyset$ in Theorem~\ref{thm::existence} simultaneously without changing Assumption~\ref{gassemilin}~\ref{gassemilin:boundedness}, since $Y_{\ad}$ defined in~\eqref{Yadptw0} and~\eqref{Yadav0} are bounded in $W^{1,q}$ again; cf.\ the discussion in Section~\ref{sec:gamma} above.
The generalization of Theorems~\ref{thm::fonintegr} and~\ref{thm::fonptw} to the additional constraints is straightforward because dealing with pointwise bounds is already possible in a less regular setting then required for dealing with gradient bounds. The FONs will contain an additional multiplier $\bar\mu \in\mathcal{M}(\bar I)$ (for~\eqref{Yadav0}) and $\bar\mu\in\mathcal{M}(\bar Q)$ (for~\eqref{Yadptw0}), respectively, and extending the proofs is straightforward; see, e.g.,~\cite{HoppeNeitzel2020} for the same quasilinear parabolic equation but without $\Fkal$. Of course, the linearized Slater conditions in Assumption~\ref{def:linearized-slater} need to replaced by ones adapted to the new $Y_{\ad}$, in particular also requiring continuity of $g$ and $g_{\avg}$.

\subsection{Pointwise bounds on each partial derivative}
The last generalization we want to discuss is the presence of bounds on each component of the gradient instead of Euclidean norms in Assumption~\ref{gasconstr}~\ref{gasconstr-ptw}. To this end,
we can modify e.g.~\eqref{Yadptw} as follows:

   \[Y_{\ad} \coloneqq \left\{ y \in \W^{1,s}(I,H^{-\zeta,p}_D, \Dkal_\zeta)\colon g_{\low} \leq \nabla y \leq g_{\up}~\text{a.e.\ on } Q \right\}\]
with $g_{\low}, g_{\up} \in L^\infty(I,L^p)^d$, $g_{\low} \leq g_{\up}$ a.e., where the inequalities have to be understood componentwise. This includes in particular the case that the Euclidean norm of the gradient in~\eqref{Yadptw} is replaced by the $\ell^\infty$-norm.

Similar considerations as above show that Theorem~\ref{thm::existence} is not affected. Let us thus only briefly explain how Theorem~\ref{thm::fonptw} has to be modified when aiming at these bounds for the gradient, of course again assuming appropriate continuity of the bounds. In this case we pick $Z = C(\overline Q)^d$, $C = \{ \varphi \in Z\colon g_{\low} \leq \varphi \leq g_{\up}\}$, and $Q = \nabla \circ \Ekal \circ S$, with $\Ekal$ the embedding $\W^{1,s}(I,H^{-\zeta,p}_D,\Dkal_\zeta) \hookrightarrow C(\overline I,C^1(\overline\Omega))$, such that $Q$ overall maps $L^s(\Lambda)$ into $C(\overline Q)^d$. Then the proof of Theorem~\ref{thm::fonptw} applies with obvious changes. The Lagrange multiplier is a vector of measures, $\bar \nu \in \Mkal(\overline Q)^d$, and the variational inequality~\eqref{fonptw_compl1} becomes equivalent to
\begin{align*}
\supp \bar \nu_i^+ \subset \bigl\{ (t,x) \in \overline Q\colon \partial_i \bar y(t,x) &= g_{\up, i}(t,x) \bigr\}, \\ \supp \bar \nu_i^- \subset \bigl\{ (t,x) \in \overline Q\colon \partial_i \bar y(t,x) & = g_{\low, i}(t,x) \bigr\},
\end{align*}
for $i=1,\dots,d$, where $\bar \nu_i = \bar \nu_i^+ - \bar \nu_i^-$ denotes the Jordan decomposition of $\bar \nu_i$. Furthermore, putting as usual
\[
\bigl\langle \nabla^* \bar \nu, \psi \bigr\rangle = \sum_{i=1}^d \langle \bar \nu_i, \partial_i \psi \rangle_{\Mkal(\overline Q), C(\overline Q)} \eqqcolon \bigl\langle -\nabla \cdot \bar\nu, \psi \bigr\rangle.
\]
for $\psi \in C(\overline I,C^1(\overline\Omega))$ and decomposing $\bar \nu = \bar \nu_Q + \delta_T \otimes \bar\nu_T$ with $\bar \nu_Q(\{T\} \times \overline\Omega) = 0$ and $\bar\nu_T \in \Mkal(\overline\Omega)^d$, the adjoint equation now reads
\[
\begin{aligned}
-\partial_t \bar w + \Akal(\bar y)^* \bar w + \Alin(\bar y)^* \bar w - \Fkal'(\bar y)^* \bar w & = \bar y - y_d - \nabla \cdot \bar \nu_Q, \\ \bar w(T) & = - \nabla \cdot \bar \nu_T.
\end{aligned}
\]

\appendix

\section{Maximal parabolic regularity}\label{sec:mpr}

In this section we collect some useful facts on maximal parabolic regularity that are used throughout the paper. We refer to~\cite{Amann1995,Amann2004,Amann2005} for several of the constructions and concepts mentioned below.

For an interval $I = (a,b) \subseteq (0,\infty)$, a number $r \in (1,\infty)$ and Banach spaces $Y \hookrightarrow_d X$, we introduce the \emph{maximal regularity space}
\[
\W^{1,r}\bigl(I, X,Y\bigr) \coloneqq W^{1,r}(I, X) \cap L^r(I, Y),
\]
equipped with the canonical sum norm; hereby, $L^r(I,Y)$ and $W^{1,r}(I,X)$ denote the underlying Bochner-Lebesgue and Bochner-Sobolev spaces.
We also use the notation $\W^{1,r}_0(I,X,Y)$ for the subset of functions in $\W^{1,r}(I,X,Y)$ that vanish at $a$, the left end point of $I$. Before we start with the actual concept of maximal parabolic regularity, we note that
\begin{equation*}
  \W^{1,r}\bigl(I,X,Y\bigr) \hookrightarrow C\bigl(\overline I,(X,Y)_{1-\frac1r,r}\bigr)
\end{equation*}
and the embedding is ``sharp'' in the sense that there is no smaller space than the real interpolation space on the right for which $t \mapsto y(t)$ is continuous for $y \in \W^{1,r}(I,X,Y)$. Thus, $(X,Y)_{1-\frac1r,r}$ is the correct function space for point evaluations in time for $\W^{1,r}(I,X,Y)$; see also Proposition~\ref{prop:cont-mpr-plus-pruess} below.

\begin{definition}[Nonautonomous maximal parabolic regularity] Let $A \colon I \to \Lkal(Y,X)$ be a family of operators on $I$ where we consider each $A(t)$ as a closed operator in $X$ with domain $Y$, such that
\begin{equation*}
  A \in L^1\bigl(I,\Lkal(Y,X)\bigr) \cap \Lkal\bigl(\W^{1,r}(I,X,Y),L^r(I,X)\bigr).
\end{equation*}
For the latter, we put $(Ay)(t) = A(t)y(t)$. Then we say that $A$ satisfies \emph{nonautonomous maximal parabolic regularity} with constant domain $Y$ on $L^r(I,X)$ iff
\begin{equation*}
  \partial_t + A \colon \W^{1,r}_0(I,X,Y) \to L^r(I,X)
\end{equation*}
is a topological isomorphism, or, equivalently, iff for every $f \in L^r(I,X)$, there is a unique solution $y \in \W^{1,r}_0(I,X,Y)$ to the abstract differential equation
\begin{equation*}
  \partial_t y(t) + A(t)y(t) = f(t)\quad \text{in}~X~\text{for}~t \in I = (a,b), \qquad y(a) = 0.
\end{equation*}
\end{definition}

We recall that it is well known that if $A \colon I \to \Lkal(Y,X)$ is in fact constant (the \emph{autonomous} case), then maximal parabolic regularity on $L^r(I,X)$ implies maximal parabolic regularity on $L^s(I,X)$ for all $s \in (1,\infty)$---see~\cite[Theorem~4.2]{Dore1993}---but that is not always true in the nonautonomous case.

Moreover, due to $I$ being a bounded interval, it is easily seen from playing with exponentials in time that a given operator $A$ satisfies maximal parabolic regularity if and only if $A + \lambda$ does so, for any scalar $\lambda$.

The crucial sufficient conditions for nonautonomous maximal parabolic regularity that we lean on in this paper are the following. We refer to~\cite[Proposition~3.1 and Theorem~5.1]{Amann2004} and the perturbation result of Pr\"uss~\cite[Corollary~3.4]{Pruess2002}:
\begin{proposition}\label{prop:cont-mpr-plus-pruess}
 Let $s \in (1,\infty)$. Suppose that
 \[A \in C\bigl(\overline I,\Lkal(Y,X)\bigr) \quad \text{and} \quad B \in L^s\bigl(I,\Lkal((X,Y)_{1-\frac1s,s},X)\bigr).\]
 Then the following are equivalent:
 \begin{enumerate}[(i)]
 \item $A$ satisfies nonautonomous maximal parabolic regularity on $L^r(I,X)$ for every $r \in (1,\infty)$ with constant domain $Y$,
 \item
 Every operator $A(\tau)$ for $\tau \in \overline I$ satisfies (autonomous) maximal parabolic regularity on $X$ with domain $Y$, \item
 The following operator is a topological isomorphism:
\begin{equation*}
  \bigl(\partial_t + A,\delta_a\bigr) \colon \W^{1,r}(I,X,Y) \to L^r(I,X) \times (X,Y)_{1-\frac1r,r}.
  \end{equation*}
 \end{enumerate}
 If any of the foregoing equivalent conditions are true, then $A+B$ also satisfies nonautonomous maximal parabolic regularity on $L^s(I,X)$.
\end{proposition}

In Proposition~\ref{prop:cont-mpr-plus-pruess}, by $\delta_a$ we mean the point evaluation $y \mapsto y(a)$ at $a$, i.e., the Dirac measure at $a$, and $a$ was the left end point of $I$.

\section{Model nonlinearity example}\label{sec:model-nonlinearity}
In this appendix we verify that our two prototypical examples for the nonlinearity $\Fkal$ stated in~\eqref{eq:examplesF1} and~\eqref{eq:examplesF2} in the introduction indeed satisfy Assumption~\ref{gassemilin}. We will demonstrate this in detail for~\eqref{eq:examplesF1} and leave an analogous discussion of~\eqref{eq:examplesF2} to the reader. We will end this appendix by pointing out some further generalizations.

\begin{lemma-and-example}\label{lem:grad-nonlinearity} Let $\zeta < 1-\frac2s$ and consider $\Fkal(y) \coloneqq \abs{\nabla y}^2$. Then $\Fkal$ satisfies Assumption~\ref{gassemilin}.
\end{lemma-and-example}

\begin{proof}
It is obvious that $\Fkal$ is a Volterra map. For the remaining assertions, let $a(v,w) \coloneqq \nabla v \cdot \nabla w$. Then $\Fkal(y) = a(y,y)$. If $q$ is such that $\frac{2d}q \leq \frac{d}p + \zeta$, then
  \begin{equation}\label{eq:quad-grad-basic-estimate}
\norm[\big]{\nabla v \cdot \nabla w}_{H^{-\zeta,p}_D} \leq C \norm[\big]{\nabla v \cdot \nabla w}_{L^{q/2}} \leq C\norm{\nabla v}_{L^q} \norm{\nabla w}_{L^q}.
  \end{equation}
  It follows that for every $r>s$, the bilinear form
   \begin{equation*} a \colon C(\overline I,W^{1,q}) \times C(\overline I,W^{1,q}) \to L^r(I,H^{-\zeta,p}_D)
  \end{equation*}
  is continuous, which in turn shows that the associated quadratic form \begin{equation*}\Fkal \colon  C(\overline I,W^{1,q}) \to L^r(I,H^{-\zeta,p}_D), \qquad \Fkal(y) \coloneqq a(y,y)\end{equation*} is continuous as well. This will imply all the assertions of Assumption~\ref{gassemilin}, provided we are able to show that for some $\frac{2d}q \leq \frac{d}p + \zeta$ as above,
\begin{equation}\label{eq:maxreg-embed-2r-gradient}
\W^{1,s}(I,H^{-\zeta,p}_D,\Dkal_\zeta) \hookrightarrow_c C(\overline I,W^{1,q}).
  \end{equation}
  To this end, we branch along the regularity cases of Assumption~\ref{assu:regularity-settings}. The general assumption for this lemma was $\zeta < 1-\frac2s$. If in addition $\zeta > \frac{d}p$, then we have $\frac{2d}p \leq \frac{d}p + \zeta$, so we can set $q \coloneqq p$ and appeal to Corollary~\ref{cor:MPR-cont-W1p-embed} to obtain~\eqref{eq:maxreg-embed-2r-gradient}. Otherwise, if $\zeta \leq \frac{d}p$, then we are in the regular case $\domLp = W^{2,p} \cap W^{1,p}_D$. Here, we choose $q>p$ close enough to $p$ such that $\frac{d}p - \frac{d}q < 1 - \frac2s - \zeta$ and in addition $\frac{2d}q \leq \frac{d}p + \zeta$. This is possible because $\frac{d}{2p} - \frac\zeta2 < 1-\frac2s - \zeta$ due to $\frac{d}p < 1- \frac2s$. Then Corollary~\ref{cor:MPR-cont-W1q-embed} strikes and yields the embedding~\eqref{eq:maxreg-embed-2r-gradient}.

  Now, having~\eqref{eq:maxreg-embed-2r-gradient} at hand, it follows directly that $\Fkal$ is well defined as a mapping $\W^{1,s}(I,H^{-\zeta,p}_D,\Dkal_\zeta) \to L^r(I,H^{-\zeta,p}_D)$ for any $r>s$ and even weak-to-strong continuous, since the embedding~\eqref{eq:maxreg-embed-2r-gradient} is compact. It is thus enough to establish the remaining assertions for $\Fkal$ defined on $C(\overline I,W^{1,q})$ which goes as follows:
  \begin{itemize}
\item Suppose that $y,z \in C(\overline I,W^{1,q})$ with their norms bounded by $M$. Then, since $a(y,y) - a(z,z) = a(y+z,y-z)$ and $a$ is continuous, \begin{equation*}\norm{\Fkal(y) - \Fkal(z)}_{L^r(I,H^{-\zeta,p}_D)} \leq 2CM \norm{y-z}_{C(\overline I,W^{1,q})},\end{equation*} so $\Fkal$ is Lipschitz continuous on bounded sets of $C(\overline I,W^{1,q})$.
\item The continuous quadratic form $\Fkal$ is Fr\'echet-differentiable with the derivative
\begin{equation}\Fkal'(y)h = 2a(y,h) = 2 \nabla y \cdot \nabla h.\label{eq:quad-grad-derivative}\end{equation}
Given $y \in C(\overline I,W^{1,q})$, we show that $\Fkal'(y)$ is an admissible perturbation for nonautonomous maximal parabolic regularity on $L^s(I,H^{-\zeta,p}_D)$ with constant domain $\Dkal_\zeta$ via Proposition~\ref{prop:cont-mpr-plus-pruess}. First, due to $(H^{-\zeta,p}_D,\Dkal_\zeta)_{1-1/s,s} \hookrightarrow [H^{-\zeta,p}_D,\Dkal_\zeta]_\theta$ for $\theta < 1-\frac1s$, it is enough to identify $\Fkal'(y)$ with an element of
$L^s(I,\Lkal([H^{-\zeta,p}_D,\Dkal_\zeta]_\theta,H^{-\zeta,p}_D))$ for such $\theta$. With the form of $\Fkal'(y)$ as in~\eqref{eq:quad-grad-derivative} and taking into account~\eqref{eq:quad-grad-basic-estimate}, this follows if $\nabla \colon [H^{-\zeta,p}_D,\Dkal_\zeta]_\theta \to L^q$ is a continuous linear operator with $q$ as considered above. The latter follows by combining Lemma~\ref{lem:complex-interpol-in-fracpower-domain} with the Kato square root property (Proposition~\ref{prop:kato}) in the irregular case $\zeta > \frac{d}p$, or, in the regular case, with Lemma~\ref{lem:regular-case-fracpower-Lipschitz} as in the proof of Corollary~\ref{cor:MPR-cont-W1q-embed} by choosing $\beta$ between $\frac{d}p - \frac{d}q$ and $1-\frac2s - \zeta$.
\item Assumption~\ref{gassemilin}~\ref{gassemilin:boundedness} is immediately true for $\Fkal$ via~\eqref{eq:quad-grad-basic-estimate} for the choice of $q$ above.\qedhere
  \end{itemize}
\end{proof}

Let us conclude with some comments on further nonlinearites that can be handled by adapting the proof of Lemma~\ref{lem:grad-nonlinearity} appropriately. First, to deal with the second prototypical nonlinearity~\eqref{eq:examplesF2} from the introduction, \[\Fkal(y)=y\beta\cdot\nabla y,\] with, say, $\beta \in L^\infty(I,L^p)^d$, one may again verify Assumption~\ref{gassemilin} by arguing quite analogously,
utilizing that $y \in \W^{1,s}(I,H^{-\zeta,p}_D,\Dkal_\zeta)$ is uniformly bounded in time and space.

 Moreover, building upon Lemma~\ref{lem:grad-nonlinearity}, we can also consider higher order nonlinearities for $\nabla y$ in the prototype form \[\Fkal_\alpha(y) \coloneqq \abs{\nabla y}^{2\alpha}\] with $\alpha > 1$ by putting $\Fkal_\alpha(y) = \abs{\Fkal(y)}^\alpha$.

 Then, in analogy to the argument based on~\eqref{eq:quad-grad-basic-estimate} in the proof of Lemma~\ref{lem:grad-nonlinearity} before, it is sufficient to determine $q$ such that \[L^{{q}/{2\alpha}} \hookrightarrow H^{-\zeta,p}_D,\] that is, $\frac{2\alpha d}{q} \leq \frac{d}p + \zeta$, and such that~\eqref{eq:maxreg-embed-2r-gradient} still holds true. We propose that this works under the assumption that \[\zeta < 1 - \frac2s - \frac{d}p + \frac{d}{p_\alpha} \quad \text{ with } \quad p_\alpha \coloneqq (2\alpha -1)p\] by choosing $q$ such that \[\frac{d}q \leq \zeta \wedge \frac{d}{p_\alpha}\quad \text{ and } \quad \frac{d}{p_\alpha} - \frac{d}q < 1 - \frac2s - \frac{d}p + \frac{d}{p_\alpha} - \zeta.\] The remaining assumptions on $\Fkal_\alpha$ to be checked then follow from the corresponding properties of $\Fkal$ and of the superposition operator $L^{q} \to L^{q/\alpha}$ induced by $\abs{\cdot}^\alpha$. We omit the details, but we point out that for $\alpha = 1$, the quantities mentioned before collapse exactly to the ones of Lemma~\ref{lem:grad-nonlinearity} and its proof.

  Finally, a similar argument as for Lemma~\ref{lem:grad-nonlinearity} allows to verify Assumption~\ref{gassemilin} for the more general variant \[\Fkal(t,y) = \eta(y) \nabla y \cdot \chi \nabla y,\] where $\eta\colon \R \rightarrow \R$ is continuously differentiable and $\chi\colon \Omega \rightarrow \R^{d \times d}$.

\section*{Acknowledgements} The authors are grateful to the anonymous referees for several comments which
have helped to improve the presentation of the paper. F.~H.\ would like to thank his former host institution University of Bonn where this research was initiated. L.~B.\ gratefully acknowledges support of his former host institution Technische Universit\"at München. This work was partially supported by DFG-Projektnummer~211504053---SFB 1060.


\medskip
Received xxxx 20xx; revised xxxx 20xx; early access xxxx 20xx.
\medskip


\begin{thebibliography}{10}

\bibitem{Amann1995} H.  Amann.
  \newblock{}{\em Linear and quasilinear parabolic problems. {V}ol. {I}}, volume~89 of {\em Monographs in
    Mathematics}.
  \newblock{}Birkh\"{a}user Boston, Inc., Boston, MA, 1995.
  \newblock{}Abstract linear theory.

\bibitem{Amann2004} H.  Amann.
  \newblock{}Maximal regularity for nonautonomous evolution equations.
  \newblock{}{\em Adv. Nonlinear Stud.}, 4(4):417--430, 2004.

\bibitem{Amann2005} H.  Amann.
  \newblock{}Nonautonomous parabolic equations involving measures.
  \newblock{}{\em Journal Math.\ Sci.} 130(4):4780--4802, 2005.

\bibitem{Amann2005_4} H.  Amann.
  \newblock{}Non-local quasi-linear parabolic equations.
  \newblock{}{\em Russ. Math. Surv.}, 60:6 1021--1033, 2005.

\bibitem{Amann2005_3} H.  Amann.
  \newblock{}Quasilinear parabolic problems via maximal regularity.
  \newblock{}{\em Adv. Differential Equations}, 10(10):1081--1110, 2005.

\bibitem{Bechtel2019} S. Bechtel and M. Egert.
  \newblock{}Interpolation theory for Sobolev functions with partially vanishing trace on irregular open
  sets.
  \newblock{}\emph{J. Fourier Anal. Appl.}, 25:2733--2781 (2019).

\bibitem{Bechtel2024} S. Bechtel.
  \newblock{}$L^p$-estimates for the square root of elliptic systems with mixed boundary conditions II\@.
  \newblock{}{\em J. Diff. Eq.} 379:104--124 (2024).

\bibitem{Bonifacius2018} L. Bonifacius and I. Neitzel.
  \newblock{}Second order optimality conditions for optimal control of quasilinear parabolic equations.
  \newblock{}{\em Math. Control Relat. Fields}, 8(1):1--34, 2018.

\bibitem{Casas86} E. Casas.
  \newblock{}Control of an elliptic problem with pointwise state constraints.
  \newblock{}{\em SIAM J. Control Optim.}, 24(6):1309--1318, 1986.

\bibitem{Casas2001} E. Casas, M. Mateos, and J.-P. Raymond.
  \newblock{}Pontryagin's principle for the control of parabolic equations with gradient state
  constraints.
  \newblock{}{\em Nonlinear Anal.}, 46(7, Ser. A\@: Theory Methods):933--956, 2001.

\bibitem{Casas1993} E. Casas.
  \newblock{}Boundary control of semilinear elliptic equations with pointwise state constraints.
  \newblock{}{\em SIAM J. Control Optim.}, 31(4):993--1006, 1993.

\bibitem{Casas2018} E. Casas and K. Chrysafinos.
  \newblock{}Analysis and optimal control of some quasilinear parabolic equations.
  \newblock{}{\em Mathematical Control and Related Fields}, 8(3--4):607--623, 2018.

\bibitem{CasasFernandez1993} E. Casas and L.A. Fern\'{a}ndez.
  \newblock{}Optimal control of semilinear elliptic equations with pointwise constraints on the gradient
  of the state.
  \newblock{}{\em Appl. Math. Optim.}, 27(1):35--56, 1993.

\bibitem{CasasFernandez1995} E. Casas and L.A. Fern\'{a}ndez.
  \newblock{}Dealing with integral state constraints in boundary control problems of quasilinear elliptic
  equations.
  \newblock{}{\em SIAM J. Control Optim.}, 33(2):568--589, 1995.

\bibitem{Casas2014} E. Casas, M. Mateos, and B. Vexler.
  \newblock{}New regularity results and improved error estimates for optimal control problems with state
  constraints.
  \newblock{}{\em ESAIM Control Optim. Calc. Var.}, 20(3):803--822, 2014.

\bibitem{Clement1992} Ph. Cl\'ement, C.\ J.\ van Duijn, Shuanhu Li.
  \newblock{}On a nonlinear elliptic-parabolic partial differential equation system in a two-dimensional
  groundwater flow problem.
  \newblock{}{\em SIAM J. Math. Anal.}, 23(4):836--851, 1992.

\bibitem{CRST93} S. Clain, J. Rappaz, M. Swierkosz and R. Touzani.
  \newblock{}Numerical modeling of induction heating for two-dimensional geometries.
  \newblock{}{\em Math. Mod. Meth. Appl. S.}, 3(6):805--822, 1993.

\bibitem{Dore1993} G. Dore.
  \newblock{}$L^p$ regularity for abstract differential equations.
  \newblock{}In: Komatsu, H.\ (eds) Functional Analysis and Related Topics, 1991. Lecture Notes in
  Mathematics, vol 1540. Springer, Berlin, Heidelberg, 1993.

\bibitem{Elschner2007} J. Elschner, J. Rehberg, and G. Schmidt.
  \newblock{}Optimal regularity for elliptic transmission problems including {$C^1$} interfaces.
  \newblock{}{\em Interfaces Free Bound.}, 9(2):233--252, 2007.

\bibitem{Disser2015} K. Disser, H-C. Kaiser and J. Rehberg.
  \newblock{}Optimal Sobolev regularity for linear second-order divergence elliptic operators occurring
  in real-world problems.
  \newblock{}{\em SIAM J. Math. Anal.} 47(3):1719--1746, 2015.

\bibitem{terElst2017} K. Disser, A.F.M. ter Elst and J. Rehberg.
  \newblock{}On maximal parabolic regularity for non-autonomous parabolic operators.
  \newblock{}{\em J. Diff. Eq.} 262(3):2039--2072 (2017).

\bibitem{Gajewski1974} H. Gajewski, K. Gr\"{o}ger, and K. Zacharias.
  \newblock{}{\em Nichtlineare {O}peratorgleichungen und {O}peratordifferentialgleichungen}.
  \newblock{}Akademie-Verlag, Berlin, 1974.
  \newblock{}Mathematische Lehrb\"{u}cher und Monographien, II.\ Abteilung, Mathematische Monographien,
  Band 38.

\bibitem{Goldberg1992} H. Goldberg, W. Kampowsky, and F. Tr\"oltzsch.
  \newblock{}On {N}emytskij operators in {$L_p$}-spaces of abstract functions.
  \newblock{}{\em Math. Nachr.}, 155:127--140, 1992.

\bibitem{Goldberg1993} H. Goldberg and F. Tr\"{o}ltzsch.
  \newblock{}Second-order sufficient optimality conditions for a class of nonlinear parabolic boundary
  control problems.
  \newblock{}{\em SIAM J. Control Optim.}, 31(4):1007--1025, 1993.

\bibitem{Grisvard1985} P. Grisvard.
  \newblock{}{\em Elliptic problems in nonsmooth domains}, volume~24 of {\em Monographs and Studies in
    Mathematics}.
  \newblock{}Pitman (Advanced Publishing Program), Boston, MA, 1985.

\bibitem{Groeger1989} K. Gr\"{o}ger.
  \newblock{}A {$W^{1,p}$}-estimate for solutions to mixed boundary value problems for second order
  elliptic differential equations.
  \newblock{}{\em Math. Ann.}, 283(4):679--687, 1989.

\bibitem{Dintelmann2009} R. Haller-Dintelmann and J. Rehberg.
  \newblock{}Maximal parabolic regularity for divergence operators including mixed boundary conditions.
  \newblock{}{\em J. Differential Equations}, 247(5):1354--1396, 2009.

\bibitem{Hesaaraki2004} M. Hesaaraki and A. Moameni.
  \newblock{}Blow-up positive solutions for a family of nonlinear parabolic equations in general domain
  in {$\mathbb{R}^N$}.
  \newblock{}{\em Michigan Math. J.}, 52(2):375--389, 2004.

\bibitem{HoppeNeitzel2020} F. Hoppe and I. Neitzel.
  \newblock{}{Optimal control of quasilinear parabolic PDEs with state constraints}.
  \newblock{}{\em SIAM J. Control Optim.}, 60(1):330--354, 2022.

\bibitem{HoppeMeinlschmidtNeitzel2023} F. Hoppe, H. Meinlschmidt, and I. Neitzel.
  \newblock{}{Global-in-time solutions and H\"older continuity for quasilinear parabolic PDEs with mixed
    boundary conditions in the Bessel dual scale}.
  \newblock{}{\em Evol. Equ. Control Theory}, early access, DOI 10.3934/eect.2024025, 2024.

\bibitem{Jonsson1984} A. Jonsson and H. Wallin.
  \newblock{}{Function Spaces on Subsets of $\mathbb{R}^n$}.
  \newblock{}Harwood Academic Publishers, 1984.

\bibitem{KPZ1986} M. Kardar, G. Parisi, and Y. Zhang.
  \newblock{} Dynamic Scaling of Growing Interfaces.
  \newblock{} {\em Phys. Rev. Lett.}. \textbf{56} (3), 889--892, 1986.

\bibitem{Kawohl1989} B. Kawohl and L.A. Peletier.
  \newblock{}Observations on blow up and dead cores for nonlinear parabolic equations.
  \newblock{}{\em Math. Z.}, 202(2):207--217, 1989.

\bibitem{Ludovici2015} F. Ludovici and W. Wollner.
  \newblock{}A priori error estimates for a finite element discretization of parabolic optimization
  problems with pointwise constraints in time on mean values of the gradient of the state.
  \newblock{}{\em SIAM J. Control Optim.}, 53(2):745--770, 2015.

\bibitem{LudoviciNeitzelWollner} F. Ludovici, I. Neitzel and W. Wollner.
  \newblock{}A priori error estimates for nonstationary optimal control problems with gradient state
  constraints.
  \newblock{}{\em PAMM}, WILEY-VCH, 15:611--612, 2015.

\bibitem{Mackenroth1983} U. Mackenroth.
  \newblock{}On parabolic distributed optimal control problems with restrictions on the gradient.
  \newblock{}{\em Appl. Math. Optim.}, 10(1):69--95, 1983.

\bibitem{Meinlschmidt2017_1} H. Meinlschmidt, C. Meyer, and J. Rehberg.
  \newblock{}Optimal control of the thermistor problem in three spatial dimensions, {P}art 1: {E}xistence
  of optimal solutions.
  \newblock{}{\em SIAM J. Control Optim.}, 55(5):2876--2904, 2017.

\bibitem{Meinlschmidt2017_2} H. Meinlschmidt, C. Meyer, and J. Rehberg.
  \newblock{}Optimal control of the thermistor problem in three spatial dimensions, {P}art 2:
  {O}ptimality conditions.
  \newblock{}{\em SIAM J. Control Optim.}, 55(4):2368--2392, 2017.

\bibitem{Meinlschmidt2016} H. Meinlschmidt and J. Rehberg.
  \newblock{}H\"older-estimates for non-autonomous parabolic problems with rough data.
  \newblock{}{\em Evol. Equ. Control Theory}, 5(1):147--184, 2016.

\bibitem{Meinlschmidt2021} H. Meinlschmidt and J. Rehberg
  \newblock{}Extrapolated elliptic regularity and application to the van Roosbroeck system of
  semiconductor equations.
  \newblock{} {\em Journal Of Differential Equations}. \textbf{280} pp. 375--404, 2021.

\bibitem{Pruess2002} J. Pr\"{u}ss.
  \newblock{}Maximal regularity for evolution equations in {$L_p$}-spaces.
  \newblock{}{\em Conf. Semin. Mat. Univ. Bari}, (285):1--39 (2003), 2002.

\end{thebibliography}
\end{document}